\documentclass[11pt,english]{article}
\usepackage[T1]{fontenc}
\usepackage[latin9]{inputenc}
\usepackage{geometry}
\usepackage{prettyref}
\usepackage{amsthm}
\usepackage{amsmath}
\usepackage{amssymb}
\usepackage{setspace}
\usepackage{esint}
\makeatletter
\theoremstyle{plain}
\newtheorem{thm}{\protect\theoremname}
  \theoremstyle{plain}
  \newtheorem{lem}[thm]{\protect\lemmaname}
  \theoremstyle{plain}
  \newtheorem{prop}[thm]{\protect\propositionname}
  \theoremstyle{plain}
  \newtheorem{cor}[thm]{\protect\corollaryname}
  \theoremstyle{remark}
  \newtheorem{rem}[thm]{\protect\remarkname}

\newrefformat{chap}{Chapter \ref{#1}}
\newrefformat{prop}{Proposition \ref{#1}}
\newrefformat{cor}{Corollary \ref{#1}}
\newrefformat{rem}{Remark \ref{#1}}
\newrefformat{ex}{Example \ref{#1}}
\newrefformat{fig}{Figure \ref{#1}}

\allowdisplaybreaks

\makeatother

\usepackage{babel}
  \providecommand{\corollaryname}{Corollary}
  \providecommand{\lemmaname}{Lemma}
  \providecommand{\propositionname}{Proposition}
  \providecommand{\remarkname}{Remark}
\providecommand{\theoremname}{Theorem}

\begin{document}

\title{Smoothness of the Orlicz Norm in\\
Musielak--Orlicz Function Spaces}

\author{Rui F.\ Vigelis%
\thanks{Computer Engineering, Campus Sobral, Federal University of Ceará,
Sobral-CE, Brazil. E-mails: rfvigelis@ufc.br.%
},\and Charles C.\ Cavalcante%
\thanks{Wireless Telecommunication Research Group, Department of Teleinformatics
Engineering, Federal University of Ceará, Fortaleza-CE, Brazil. E-mail:
charles@ufc.br.%
}}
\maketitle
\begin{abstract}
In this paper, we present a characterization of support functionals
and smooth points in $L_{0}^{\Phi}$, the Musielak--Orlicz space equipped
with the Orlicz norm. As a result, criterion for the smoothness of
$L_{0}^{\Phi}$ is also obtained. Some expressions involving the norms
of functionals in $(L_{0}^{\Phi})^{*}$, the topological dual of $L_{0}^{\Phi}$,
are proved for arbitrary Musielak--Orlicz functions.\textit{}
\end{abstract}

\section{Introduction}

Characterization of support functionals and smooth points, as well
as criterion for the smoothness of Musielak--Orlicz (function) spaces
equipped with the Orlicz norm, which we denote by $L_{0}^{\Phi}$,
are already known \cite{Hudzik:1997} when the Musielak--Orlicz function
$\Phi$ is finite-valued and $\Phi(t,u)/u\rightarrow\infty$ as $u\rightarrow\infty$
for $\mu$-a.e.\ $t\in T$. In this paper, we show these results
for arbitrary Musielak--Orlicz functions, which can take values in
the extended real numbers. The proofs follow the main lines of the
paper \cite{Hudzik:1997}, with improvements. For instance, we have
neither used the Bishop--Phelps Theorem \cite{Phelps:1993}, nor the
concept of measurable selectors \cite{Aliprantis:2006a}. (As a consequence,
see the functions $u_{*}$ and $u^{*}$ constructed in \prettyref{rem:u_stars}.)
To find these characterizations, some expressions involving the norms
of functionals in $(L_{0}^{\Phi})^{*}$, the topological dual of $L_{0}^{\Phi}$,
are extended to arbitrary Musielak--Orlicz functions. In the proof
of these extensions, a strategy consisted of taking a sequence of
finite-valued Musielak--Orlicz functions converging upward to an arbitrary
Musielak--Orlicz function (see \prettyref{lem:seqMO}). Next we present
some definitions and results related to Musielak--Orlicz spaces \cite{Rao:1991,Musielak:1983}.

Let $(T,\Sigma,\mu)$ be a non-atomic, $\sigma$-finite measure space.
We say that $\Phi\colon T\times[0,\infty]\rightarrow[0,\infty]$ is
a \textit{Musielak--Orlicz function} if, for $\mu$-a.e.\ $t\in T$,
\begin{itemize}
\item [(i)] $\Phi(t,\cdot)$ is convex and lower semi-continuous,
\item [(ii)] $\Phi(t,0)=\lim_{u\downarrow0}\Phi(t,u)=0$ and $\Phi(t,\infty)=\infty$,
\item [(iii)] $\Phi(\cdot,u)$ is measurable for all $u\geq0$.
\end{itemize}
The \textit{complementary function} $\Phi^{*}\colon T\times[0,\infty]\rightarrow[0,\infty]$
to a Musielak--Orlicz function $\Phi$ is defined by
\begin{equation}
\Phi^{*}(t,v)=\sup_{u>0}(uv-\Phi(t,u)),\quad\text{for all }v\geq0.\label{eq:fenchel-s}
\end{equation}
It can be verified that $\Phi^{*}$ is a Musielak--Orlicz function.
Given any Musielak--Orlicz function $\Phi$, we denote $\partial\Phi(t,u)=[\Phi_{-}'(t,u),\Phi_{+}'(t,u)]$,
where $\Phi_{-}'(t,u)$ and $\Phi_{+}'(t,u)$ are the left- and right-derivatives
of $\Phi(t,\cdot)$ at any $u\geq0$. The functions $\Phi$ and $\Phi^{*}$
satisfy the \textit{Young's inequality}
\begin{equation}
uv\leq\Phi(t,u)+\Phi^{*}(t,v),\quad\text{for all }u,v\geq0,\label{eq:YoungIneq}
\end{equation}
which reduces to an equality when $v\in\partial\Phi(t,u)$ if $u$
is given, or when $u\in\partial\Phi^{*}(t,v)$ if $v$ is given. We
also define the functions $a_{\Phi}(t)=\sup\{u\geq0:\Phi(t,u)=0\}$
and $b_{\Phi}(t)=\sup\{u\geq0:\Phi(t,u)<\infty\}$.

Let $L^{0}$ denote the space of all real-valued measurable functions
on $T$, with equality $\mu$-a.e. Given a Musielak--Orlicz function
$\Phi$, we define the functional
\begin{equation}
I_{\Phi}(u)=\int_{T}\Phi(t,|u(t)|)d\mu,\quad\text{for any }u\in L^{0}.\label{eq:I_Phi}
\end{equation}
The \textit{Musielak--Orlicz (function) space}, \textit{Musielak--Orlicz
(function) class}, and \textit{(function) space of finite elements}
are given by 
\begin{align*}
L^{\Phi} & =\{u\in L^{0}:I_{\Phi}(\lambda u)<\infty\text{ for some }\lambda>0\},\\
\tilde{L}^{\Phi} & =\{u\in L^{0}:I_{\Phi}(u)<\infty\},\\
\intertext{and}E^{\Phi} & =\{u\in L^{0}:I_{\Phi}(\lambda u)<\infty\text{ for all }\lambda>0\},
\end{align*}
respectively. Clearly, $E^{\Phi}\subseteq\tilde{L}^{\Phi}\subseteq L^{\Phi}$.
The Musielak--Orlicz space $L^{\Phi}$ can be viewed as the smallest
vector subspace of $L^{0}$ that contains $\tilde{L}^{\Phi}$, and
$E^{\Phi}$ is the largest vector subspace of $L^{0}$ that is contained
in $\tilde{L}^{\Phi}$.

The Musielak--Orlicz space $L^{\Phi}$ is a Banach space when it is
endowed with any of the norms: 
\begin{align}
\Vert u\Vert_{\Phi} & =\inf\Bigl\{\lambda>0:I_{\Phi}\Bigl(\frac{u}{\lambda}\Bigr)\leq1\Bigr\},\label{eq:norm_Luxemburg}\\
\Vert u\Vert_{\Phi,0} & =\sup\biggl\{\biggl|\int_{T}uvd\mu\biggr|:v\in\tilde{L}^{\Phi^{*}}\text{ and }I_{\Phi^{*}}(v)\leq1\biggr\},\label{eq:norm_Orlicz}\\
\intertext{and}\Vert u\Vert_{\Phi,A} & =\inf_{k>0}\frac{1}{k}(1+I_{\Phi}(ku)),\label{eq:norm_Amemiya}
\end{align}
which are called the \textit{Luxemburg}, \textit{Orlicz} and \textit{Amemiya
norms}, respectively. The Musielak--Orlicz space equipped with the
Orlicz norm is denoted by $L_{0}^{\Phi}$. The Luxemburg and Orlicz
norms are equivalent and are related by the inequalities $\Vert u\Vert_{\Phi}\leq\Vert u\Vert_{\Phi,0}\leq2\Vert u\Vert_{\Phi}$,
for any $u\in L^{\Phi}$. In addition, as shown in \cite{Hudzik:2000,Fan:2007},
the Orlicz and Amemiya norms coincides, i.e., $\Vert u\Vert_{\Phi,0}=\Vert u\Vert_{\Phi,A}$,
for all $u\in L^{\Phi}$. The Amemiya norm is a special case of the
$p$-Amemiya norm for $p=1$. For more details on $p$-Amemiya norms
we refer to \cite{Cui:2008}. 

For any $u\in L^{\Phi}$, we denote by $K(u)$ the set of all $k>0$
for which the infimum in \prettyref{eq:norm_Amemiya} is attained.
If $I_{\Phi^{*}}(b_{\Phi^{*}}\chi_{\operatorname{supp}u})>1$, where
$\operatorname{supp}u=\{t\in T:|u(t)|>0\}$, and we denote
\begin{align*}
k_{u}^{*} & =\inf\{k>0:I_{\Phi^{*}}(\Phi_{+}'(t,|ku(t)|))\geq1\},\\
k_{u}^{**} & =\sup\{k>0:I_{\Phi^{*}}(\Phi_{+}'(t,|ku(t)|))\leq1\},
\end{align*}
then $0<k_{u}^{*}\leq k_{u}^{**}<\infty$, and $\Vert u\Vert_{\Phi,0}=\frac{1}{k}(1+I_{\Phi}(ku))$
if and only if $k\in[k_{u}^{*},k_{u}^{**}]$. If $I_{\Phi^{*}}(b_{\Phi^{*}}\chi_{\operatorname{supp}u})\leq1$,
then $\Vert u\Vert_{\Phi,0}=\int_{T}|u|b_{\Phi^{*}}d\mu$. 

If we can find a non-negative function $f\in\tilde{L}^{\Phi}$ and
a constant $K>0$ such that
\begin{equation}
\Phi(t,2u)\leq K\Phi(t,u),\quad\text{for all }u\geq f(t),\label{eq:D2_0}
\end{equation}
then we say that $\Phi$ satisfies the \textit{$\Delta_{2}$-condition},
or belongs to the \textit{$\Delta_{2}$-class} (denoted by $\Phi\in\Delta_{2}$).
The spaces $E^{\Phi}$ and $L^{\Phi}$ coincide when $\Phi$ satisfies
the $\Delta_{2}$-condition. On the other hand, if $\Phi$ is finite-valued
and does not satisfy the $\Delta_{2}$-condition, then the Musielak--Orlicz
class $\tilde{L}^{\Phi}$ is not open and its interior coincides with
\[
B_{0}(E^{\Phi},1)=\{u\in L^{\Phi}:\inf_{v\in E^{\Phi}}\Vert u-v\Vert_{\Phi,0}<1\},
\]
or, equivalently, $B_{0}(E^{\Phi},1)\varsubsetneq\tilde{L}^{\Phi}\varsubsetneq\overline{B}_{0}(E^{\Phi},1)$.

Musielak--Orlicz spaces are endowed with the structure of Banach lattices
\cite{Aliprantis:2006}. This property can be used in a more refined
analysis of the (topological) dual space of $L^{\Phi}$, which is
denoted by $(L^{\Phi})^{*}$. Fixed $v\in L^{\Phi^{*}}$, the expression
\begin{equation}
f_{v}(u):=\int_{T}uvd\mu,\qquad\text{for all }u\in L^{\Phi},\label{eq:f_v}
\end{equation}
defines a functional in $(L^{\Phi})^{*}$. A functional that can be
written in the form \prettyref{eq:f_v} is said to be \textit{order
continuous}. Unless the Musielak--Orlicz function $\Phi$ satisfies
the $\Delta_{2}$-condition, not all functionals in $(L^{\Phi})^{*}$
are represented by \prettyref{eq:f_v} for some $v\in L^{\Phi^{*}}$.
However, every functional $f\in(L^{\Phi})^{*}$ can be uniquely expressed
as 
\[
f=f_{\mathrm{c}}+f_{\mathrm{s}},
\]
where $f_{\mathrm{c}}$ and $f_{\mathrm{s}}$ are said to be the \textit{order
continuous} (i.e., $f_{\mathrm{c}}=f_{v}$ for some $v\in L^{\Phi^{*}}$)
and \textit{singular component} of $f$, respectively. A functional
$f\in(L^{\Phi})^{*}$ for which $f_{\mathrm{c}}=0$ is called \textit{purely
singular}. If the Musielak--Orlicz function $\Phi$ is finite-valued,
then purely singular functionals are characterized as those functionals
vanishing on $E^{\Phi}$. Unfortunately, this characterization can
not be used if the Musielak--Orlicz function $\Phi$ is not finite-valued,
since we can have $E^{\Phi}=\{0\}$. We say that a functional $f\in(L^{\Phi})^{*}$
is \textit{positive} if $f(u)\geq0$ for all non-negative functions
$u\in L^{\Phi}$. To find the order continuous and singular component
of a positive functional $f\in(L^{\Phi})^{*}$, we can use 
\begin{align}
f_{\mathrm{c}}(u) & =\inf\{\sup_{n\geq1}f(u_{n}):0\leq u_{n}\uparrow u\}\label{eq:f_c}\\
\intertext{and}f_{\mathrm{s}}(u) & =\sup\{\inf_{n\geq1}f(u_{n}):u\geq u_{n}\downarrow0\},\label{eq:f_s}
\end{align}
for any non-negative functions $u\in L^{\Phi}$. Expressions \prettyref{eq:f_c}
and \prettyref{eq:f_s} are valid for arbitrary Musielak--Orlicz functions.
For any $f\in(L^{\Phi})^{*}$, we define the norms 
\[
\Vert f\Vert_{0}=\sup_{u\in L^{\Phi}}\frac{|f(u)|}{\Vert u\Vert_{\Phi}},\qquad\text{and}\qquad\Vert f\Vert=\sup_{u\in L_{0}^{\Phi}}\frac{|f(u)|}{\Vert u\Vert_{\Phi,0}}.
\]
Thanks to \prettyref{eq:f_c} and \prettyref{eq:f_s}, we can show,
in \prettyref{sec:auxiliary_results}, some results related to the
norms of functionals in $(L^{\Phi})^{*}$ for arbitrary Musielak--Orlicz
functions. Assuming that the Musielak--Orlicz function $\Phi$ is
finite-valued, one can verify that $(E^{\Phi})^{*}\simeq L^{\Phi^{*}}$. 

Let $(X,\Vert\cdot\Vert)$ be a Banach space, whose (topological)
dual space is denoted by $X^{*}$. A \textit{support functional} at
$x\in X\setminus\{0\}$ is a norm-one functional $f\in X^{*}$ such
that $f(x)=\Vert x\Vert$. The Hahn--Banach Theorem ensures the existence
of at least one support functional. If there exists only one support
functional at $x\in X\setminus\{0\}$, then $x$ is said to be a \textit{smooth
point}. A Banach space $X$ is called \textit{smooth} if every $x\in X\setminus\{0\}$
is a smooth point. 

The rest of this paper is organized as follows. In \prettyref{sec:auxiliary_results},
some results related to the norms of functionals in $(L^{\Phi})^{*}$
are proved for arbitrary Musielak--Orlicz functions. In \prettyref{sec:main_results},
characterization of support functionals and smooth points in $L_{0}^{\Phi}$
are presented, and we obtain necessary and sufficient conditions for
the smoothness of $L_{0}^{\Phi}$.

\section{Auxiliary results\label{sec:auxiliary_results}}

In this section, some expressions involving the norms of functionals
in $(L^{\Phi})^{*}$ are proved for arbitrary Musielak--Orlicz functions.
To show these results, we make use of the lemmas below.

\begin{lem}
\label{lem:seq_un_IPhiS} Let $\Phi$ be an arbitrary Musielak--Orlicz
function. If $u\colon T\rightarrow[0,\infty)$ is a measurable function
satisfying $\Phi(t,u(t))<\infty$ for $\mu$-a.e.\ $t\in T$, then
we can find a sequence of non-negative measurable functions $\{u_{n}\}$
such that $u_{n}\uparrow u$, and $I_{\Phi^{*}}(\Phi_{-}'(t,u_{n}(t)))<\infty$
for all $n\geq1$.\end{lem}
\begin{proof}
For each $n\geq1$, define $\widetilde{u}_{n}(t)=\max(0,u(t)-1/n)$.
In view of $\widetilde{u}_{n}<b_{\Phi}$, it follows that $\Phi^{*}(t,\Phi_{-}'(t,\widetilde{u}_{n}(t)))<\infty$
for $\mu$-a.e.\ $t\in T$. Let $\{T_{n}\}$ be a non-decreasing
sequence of measurable sets such that $0<\mu(T_{n})<\infty$ and $\mu(T\setminus\bigcup_{n=1}^{\infty}T_{n})=0$.
We can find, for each $n\geq1$, a sufficiently large $m_{n}\geq1$
such that the set $A_{n}=\{t\in T_{n}:\Phi^{*}(t,\Phi_{-}'(t,\widetilde{u}_{n}(t)))\leq m_{n}\}$
satisfies $\mu(T_{n}\setminus A_{n})<2^{-n}$. Let $B_{n}=\bigcap_{k=n}^{\infty}A_{k}$.
Clearly, $B_{n}\subseteq B_{n+1}$ and $B_{n}\subseteq T_{n}$ for
all $n\geq1$. Thus, for any $n\geq m$, we can write 
\begin{align*}
\mu(T_{m}\setminus B_{n}) & =\mu\biggl(\bigcup_{k=n}^{\infty}T_{m}\setminus A_{k}\biggr)\leq\sum_{k=n}^{\infty}\mu(T_{m}\setminus A_{k})\\
 & \leq\sum_{k=n}^{\infty}\mu(T_{k}\setminus A_{k})\leq\sum_{k=n}^{\infty}2^{-k}\\
 & =2^{-n+1},
\end{align*}
from which we can conclude that $\mu(T\setminus\bigcup_{n=1}^{\infty}B_{n})=0$.
Defining $u_{n}=\widetilde{u}_{n}\chi_{B_{n}}$, we obtain that $u_{n}\uparrow u$,
and $I_{\Phi^{*}}(\Phi_{-}'(t,u_{n}(t)))\leq m_{n}\mu(T_{n})<\infty$
for all $n\geq1$.
\end{proof}

\begin{lem}
\label{lem:seqMO} Let $\Phi$ be an arbitrary Musielak--Orlicz function.
\begin{itemize}
\item [\normalfont(a)] Then there exists a non-decreasing sequence of finite-valued
Musielak--Orlicz functions $\{\Phi_{n}\}$ converging upward to $\Phi$,
i.e., such that $\Phi_{n}(t,u)\uparrow\Phi(t,u)$, for all $u\geq0$,
and $\mu$-a.e.\ $t\in T$. 
\item [\normalfont(b)] In addition, for any such sequence $\{\Phi_{n}\}$,
if $u$ is a function belonging to $L^{\Phi_{n}}$ for all $n\geq1$,
and the sequence $\{\Vert u\Vert_{\Phi_{n}}\}$ is bounded, then $u$
belongs to $L^{\Phi}$, and $\Vert u\Vert_{\Phi_{n}}\uparrow\Vert u\Vert_{\Phi}$.
\end{itemize}
\end{lem}
\begin{proof}
(a) For each $n\geq1$, we define the Musielak--Orlicz function 
\[
\Phi_{n}(t,u)=\int_{0}^{u}\min(\Phi_{-}'(t,x),n)dx.
\]
Clearly, $\Phi_{n}(t,u)=\Phi(t,u)$ for any $u\geq0$ satisfying $\Phi_{-}'(t,u)\leq n$.
In addition, $\Phi_{n}(t,u)\uparrow\infty$ for any $u>0$ such that
$\Phi_{-}'(t,u)=\infty$. Thus $\Phi_{n}(t,u)\uparrow\Phi(t,u)$ for
all $u\geq0$.

(b) The case $u=0$ is trivial. So we assume that $u\neq0$. Since
$I_{\Phi_{m}}(u/\lambda)\leq I_{\Phi_{n}}(u/\lambda)$ for any $\lambda>0$
and $m\leq n$, it follows that $\Vert u\Vert_{\Phi_{m}}\leq\Vert u\Vert_{\Phi_{n}}$.
Thus the sequence $\{\Vert u\Vert_{\Phi_{n}}\}$ converges upward
to some $c>0$. In view of Fatou's Lemma, for any $\lambda>c$, we
have that 
\[
I_{\Phi}(u/\lambda)\leq\liminf_{n\rightarrow\infty}I_{\Phi_{n}}(u/\lambda)\leq1.
\]
Hence $u\in L^{\Phi}$ and $\Vert u\Vert_{\Phi}\leq c$. Now, for
any $\lambda<c$, and a sufficiently large $n\geq1$ such that $\Vert u\Vert_{\Phi_{n}}>\lambda$,
we obtain that $I_{\Phi}(u/\lambda)\geq I_{\Phi_{n}}(u/\lambda)>1$.
Consequently, $\Vert u\Vert_{\Phi}=c$.
\end{proof}

\begin{prop}
\label{prop:OrliczLuxemExpr} The Orlicz and Luxemburg norms can be
expressed respectively as
\begin{align}
\Vert u\Vert_{\Phi,0} & =\sup\biggl\{\biggl|\int_{T}uvd\mu\biggr|:v\in L^{\Phi^{*}}\text{ and }\Vert v\Vert_{\Phi^{*}}\leq1\biggr\}\label{eq:OrliczNorm1}\\
\intertext{and}\Vert u\Vert_{\Phi} & =\sup\biggl\{\biggl|\int_{T}uvd\mu\biggr|:v\in L^{\Phi^{*}}\text{ and }\Vert v\Vert_{\Phi^{*},0}\leq1\biggr\}.\label{eq:LuxemburgNorm1}
\end{align}
\end{prop}
\begin{proof}
By the definition of Luxemburg norm, it is clear that $v\in L^{\Phi^{*}}$
satisfies $\Vert v\Vert_{\Phi^{*}}\leq1$ if and only if $I_{\Phi^{*}}(v)\leq1$.
Therefore, the Orlicz norm can be expressed as in \prettyref{eq:OrliczNorm1}.
A consequence of \prettyref{eq:OrliczNorm1} is Hölder's Inequality:
\[
\biggl|\int_{T}uvd\mu\biggr|\leq\Vert u\Vert_{\Phi}\Vert v\Vert_{\Phi^{*},0},
\]
which is employed in the proof of \prettyref{eq:LuxemburgNorm1}.

We will show that \prettyref{eq:LuxemburgNorm1} holds. First, we
assume that $\Phi$ is finite-valued. Without loss of generality,
we also assume that $u\geq0$ and the supremum in \prettyref{eq:LuxemburgNorm1}
is equal to $1$. From Hölder's Inequality, it follows that 
\begin{equation}
1=\sup\biggl\{\biggl|\int_{T}uvd\mu\biggr|:v\in L^{\Phi^{*}}\text{ and }\Vert v\Vert_{\Phi^{*},0}\leq1\biggr\}\leq\Vert u\Vert_{\Phi}.\label{eq:OrliczNormRestric_ineq}
\end{equation}
Suppose that the last inequality in \prettyref{eq:OrliczNormRestric_ineq}
is strict. According to \prettyref{lem:seq_un_IPhiS}, there exists
a sequence of non-negative measurable functions $\{u_{n}\}$ such
that $u_{n}\uparrow u$ and $I_{\Phi^{*}}(\Phi_{-}'(t,u_{n}(t)))<\infty$
for each $n\geq1$. Since $\Vert u\Vert_{\Phi}>1$, we obtain that
$I_{\Phi}(u)>1$. Then we can find a sufficiently large $n_{0}\geq1$
for which the function $u_{0}:=u_{n_{0}}$ satisfies $I_{\Phi}(u_{0})>1$.
Define the function
\[
v_{0}(t)=\frac{\Phi_{-}'(t,u_{0}(t))}{1+I_{\Phi^{*}}(\Phi_{-}'(t,u_{0}(t)))},
\]
which belongs to $\tilde{L}^{\Phi^{*}}$. For any non-negative function
$w\in L^{\Phi}$ such that $I_{\Phi}(w)\leq1$, it follows that 
\[
\int_{T}wv_{0}d\mu\leq\frac{I_{\Phi}(w)+I_{\Phi^{*}}(\Phi_{-}'(t,u_{0}(t)))}{1+I_{\Phi^{*}}(\Phi_{-}'(t,u_{0}(t)))}\leq1.
\]
Hence $\Vert v_{0}\Vert_{\Phi^{*},0}\leq1$. In addition, we can write
\begin{align*}
\int_{T}uv_{0}d\mu & \geq\int_{T}u_{0}v_{0}d\mu\\
 & =\frac{I_{\Phi}(u_{0})+I_{\Phi^{*}}(\Phi_{-}'(t,u_{0}(t)))}{1+I_{\Phi^{*}}(\Phi_{-}'(t,u_{0}(t)))}\\
 & >1,
\end{align*}
which is a contradiction to \prettyref{eq:OrliczNormRestric_ineq}.
Therefore, the last inequality in \prettyref{eq:OrliczNormRestric_ineq}
cannot be strict.

Now assume that $\Phi$ is arbitrary. According to \prettyref{lem:seqMO}--(a),
we can find a non-decreasing sequence of finite-valued Musielak--Orlicz
functions $\{\Phi_{n}\}$ converging upward to $\Phi$. The inequality
$\Phi^{*}(t,v)\leq\Phi_{n}^{*}(t,v)$, for all $v\geq0$ and $\mu$-a.e.\ $t\in T$,
implies that $L^{\Phi_{n}^{*}}\subseteq L^{\Phi^{*}}$. Moreover,
for any $v\in L^{\Phi_{n}^{*}}$, 
\begin{align*}
\Vert v\Vert_{\Phi^{*},0} & =\sup\biggl\{\biggl|\int_{T}uvd\mu\biggr|:u\in\tilde{L}^{\Phi}\text{ and }I_{\Phi}(u)\leq1\biggr\}\\
 & \leq\sup\biggl\{\biggl|\int_{T}uvd\mu\biggr|:u\in\tilde{L}^{\Phi_{n}}\text{ and }I_{\Phi_{n}}(u)\leq1\biggr\}\\
 & =\Vert v\Vert_{\Phi_{n}^{*},0}.
\end{align*}
Consequently, we can write
\begin{align*}
\Vert u\Vert_{\Phi_{n}} & =\sup\biggl\{\biggl|\int_{T}uvd\mu\biggr|:v\in L^{\Phi_{n}^{*}}\text{ and }\Vert v\Vert_{\Phi_{n}^{*},0}\leq1\biggr\}\\
 & \leq\sup\biggl\{\biggl|\int_{T}uvd\mu\biggr|:v\in L^{\Phi^{*}}\text{ and }\Vert v\Vert_{\Phi^{*},0}\leq1\biggr\}\\
 & \leq\Vert u\Vert_{\Phi}.
\end{align*}
In view of \prettyref{lem:seqMO}--(b), the convergence $\Vert u\Vert_{\Phi_{n}}\uparrow\Vert u\Vert_{\Phi}$
implies expression \prettyref{eq:LuxemburgNorm1}.
\end{proof}

For any $u\in L^{\Phi}$, we define 
\begin{align*}
\theta_{\Phi}(u) & =\inf\{\lambda>0:I_{\Phi}(u/\lambda)<\infty\}\\
\intertext{and}Q_{\Phi}(u) & =\sup\{\inf_{n\geq1}\Vert u_{n}\Vert_{\Phi}:|u|\geq u_{n}\downarrow0\},\\
Q_{\Phi,0}(u) & =\sup\{\inf_{n\geq1}\Vert u_{n}\Vert_{\Phi,0}:|u|\geq u_{n}\downarrow0\}.
\end{align*}
These functionals are related to the norms of purely singular functionals.
A remarkable property is that these functionals coincide. To show
this claim, we need the following lemma.

\begin{lem}
\label{lem:IPhi_un_infty} Let $u\in L^{\Phi}$ be such that $I_{\Phi}(u)=\infty$.
Then there exists a sequence of measurable functions $\{u_{n}\}$
such that $|u|\geq u_{n}\downarrow0$ and $I_{\Phi}(u_{n})=\infty$
for all $n\geq1$.\end{lem}
\begin{proof}
Let $B=\{t\in T:\Phi(t,|u(t)|)=\infty\}$. First we assume that $\mu(B)>0$.
Then we can find a non-decreasing sequence of measurable sets $\{B_{n}\}$,
with positive measure $\mu(B_{n})>0$, and such that $B_{n}\subseteq B$
and $\mu(B_{n})\downarrow0$. For each $n\geq1$, define the functions
$u_{n}=|u|\chi_{B_{n}}$. Clearly, $|u|\geq u_{n}\downarrow0$ and
$I_{\Phi}(u_{n})=\infty$ for all $n\geq1$. Now suppose that $\mu(B)=0$.
Let $\{T_{n}\}$ be a non-decreasing sequence of measurable sets such
that $0<\mu(T_{n})<\infty$ and $\mu(T\setminus\bigcup_{n=1}^{\infty}T_{n})=0$.
Define $A_{n}=\{t\in T_{n}:\Phi(t,|u(t)|)\leq n\}$. Clearly, the
sequence $\{A_{n}\}$ is non-decreasing, and $\mu(T\setminus\bigcup_{n=1}^{\infty}A_{n})=0$.
For each $n\geq1$, we define the functions $u_{n}=|u|\chi_{T\setminus A_{n}}$,
which satisfy $|u|\geq u_{n}\downarrow0$. Observing that 
\[
\infty=I_{\Phi}(u)=I_{\Phi}(u_{n})+I_{\Phi}(u\chi_{A_{n}})\leq I_{\Phi}(u_{n})+n\mu(T_{n}),
\]
we conclude that $I_{\Phi}(u_{n})=\infty$ for all $n\geq1$.
\end{proof}

\begin{prop}
\label{prop:theta_Q_Q0} For every $u\in L^{\Phi}$, there holds that
$\theta_{\Phi}(u)=Q_{\Phi}(u)=Q_{\Phi,0}(u)$.\end{prop}
\begin{proof}
It is clear that $Q_{\Phi}(u)\leq Q_{\Phi,0}(u)$. Fix any $\varepsilon>0$.
Let $\{u_{n}\}$ be a sequence in $L^{\Phi}$ such that $|u|\geq u_{n}\downarrow0$
and $Q_{\Phi,0}(u)-\varepsilon\leq\inf_{n\geq1}\Vert u_{n}\Vert_{\Phi,0}$.
Take any $\lambda>\theta_{\Phi}(u)$. In view of $I_{\Phi}(u/\lambda)<\infty$,
we obtain that $I_{\Phi}(u_{n}/\lambda)\downarrow0$. Hence 
\[
Q_{\Phi,0}(u)-\varepsilon\leq\inf_{n\geq1}\Vert u_{n}\Vert_{\Phi,0}\leq\inf_{n\geq1}\lambda(1+I_{\Phi}(u_{n}/\lambda))=\lambda.
\]
Since $\varepsilon>0$ and $\lambda>\theta_{\Phi}(u)$ are arbitrary,
we have that $Q_{\Phi,0}(u)\leq\theta_{\Phi}(u)$. If $\theta_{\Phi}(u)=0$
then $Q_{\Phi}(u)=Q_{\Phi,0}(u)=\theta_{\Phi}(u)=0$. So we assume
that $\theta_{\Phi}(u)>0$. Now let $\lambda<\theta_{\Phi}(u)$. Clearly,
$I_{\Phi}(u/\lambda)=\infty$. According to \prettyref{lem:IPhi_un_infty},
we can find a sequence of measurable functions $\{u_{n}\}$ such that
$|u|\geq u_{n}\downarrow0$ and $I_{\Phi}(u_{n}/\lambda)=\infty$
for all $n\geq1$. From the definition of Luxemburg norm, it follows
that $\Vert u_{n}\Vert_{\Phi}\geq\lambda$ for all $n\geq1$. Then
we can write $Q_{\Phi}(u)\geq\inf_{n\geq1}\Vert u_{n}\Vert_{\Phi}\geq\lambda$.
Because $\lambda<\theta_{\Phi}(u)$ is arbitrary, we conclude that
$Q_{\Phi}(u)\geq\theta_{\Phi}(u)$. Therefore, $\theta_{\Phi}(u)=Q_{\Phi}(u)=Q_{\Phi,0}(u)$.
\end{proof}

\begin{prop}
\label{prop:norm_sing_func} If the functional $f\in(L^{\Phi})^{*}$
is purely singular, then
\begin{equation}
\Vert f\Vert_{0}=\sup_{u\in L^{\Phi}}\frac{|f(u)|}{Q_{\Phi}(u)},\qquad\text{and}\qquad\Vert f\Vert=\sup_{u\in L^{\Phi}}\frac{|f(u)|}{Q_{\Phi,0}(u)},\label{eq:norm_sing_func_1}
\end{equation}
or, equivalently, 
\begin{equation}
\Vert f\Vert_{0}=\Vert f\Vert=\sup_{u\in\tilde{L}^{\Phi}}|f(u)|=\sup_{u\in L^{\Phi}}\frac{|f(u)|}{\theta_{\Phi}(u)}.\label{eq:norm_sing_func_2}
\end{equation}
\end{prop}
\begin{proof}
Without loss of generality, we can assume that $f\geq0$. The equivalence
between \prettyref{eq:norm_sing_func_1} and \prettyref{eq:norm_sing_func_2}
follows from \prettyref{prop:theta_Q_Q0}. Since $u\in\tilde{L}^{\Phi}$
if $\Vert u\Vert_{\Phi}\leq1$, and $\theta_{\Phi}(u)\leq1$ for any
$u\in\tilde{L}^{\Phi}$, we can write
\begin{align}
\Vert f\Vert_{0} & =\sup_{u\in L^{\Phi}}\frac{|f(u)|}{\Vert u\Vert_{\Phi}}\leq\sup_{u\in\tilde{L}^{\Phi}}|f(u)|\nonumber \\
 & \leq\sup_{u\in\tilde{L}^{\Phi}}\frac{|f(u)|}{\theta_{\Phi}(u)}\leq\sup_{u\in L^{\Phi}}\frac{|f(u)|}{\theta_{\Phi}(u)}\nonumber \\
 & =\sup_{u\in L_{+}^{\Phi}}\frac{f(u)}{\theta_{\Phi}(u)}.\label{eq:futhetaPhiu}
\end{align}
Now, for any $u\in L_{+}^{\Phi}$, we have that 
\begin{align}
f(u) & =\sup\{\inf\nolimits _{n\geq1}f(u_{n}):u\geq u_{n}\downarrow0\}\nonumber \\
 & \leq\sup\{\inf\nolimits _{n\geq1}\Vert f\Vert\Vert u_{n}\Vert_{\Phi,0}:u\geq u_{n}\downarrow0\}\nonumber \\
 & =\Vert f\Vert\sup\{\inf\nolimits _{n\geq1}\Vert u_{n}\Vert_{\Phi,0}:u\geq u_{n}\downarrow0\}\nonumber \\
 & =\Vert f\Vert Q_{\Phi,0}(u)=\Vert f\Vert\theta_{\Phi}(u).\label{eq:norm_sing_func}
\end{align}
From \prettyref{eq:futhetaPhiu} and \prettyref{eq:norm_sing_func},
it follows that $\Vert f\Vert_{0}\leq\Vert f\Vert$. Because $\Vert f\Vert\leq\Vert f\Vert_{0}$
is also satisfied, we obtain \prettyref{eq:norm_sing_func_2}.
\end{proof}

\begin{prop}
Every functional $f=f_{\mathrm{c}}+f_{\mathrm{s}}\in(L^{\Phi})^{*}$
satisfies $\Vert f\Vert_{0}=\Vert f_{\mathrm{c}}\Vert_{0}+\Vert f_{\mathrm{s}}\Vert_{0}$.\end{prop}
\begin{proof}
Since $|f|_{\mathrm{c}}=|f_{\mathrm{c}}|$ and $|f|_{\mathrm{s}}=|f_{\mathrm{s}}|$,
we can assume that $f\geq0$. Clearly, $\Vert f\Vert_{0}\leq\Vert f_{\mathrm{c}}\Vert_{0}+\Vert f_{\mathrm{s}}\Vert_{0}$.
Given any $\varepsilon>0$, we select non-negative functions $u,v\in L^{\Phi}$
with $\Vert u\Vert_{\Phi}\leq1$ and $\Vert v\Vert_{\Phi}\leq1$ such
that 
\[
f_{\mathrm{c}}(u)\geq\Vert f_{\mathrm{c}}\Vert_{0}-\varepsilon,\qquad\text{and}\qquad f_{\mathrm{s}}(v)\geq\Vert f_{\mathrm{s}}\Vert_{0}-\varepsilon.
\]
In view of \prettyref{eq:f_s}, there exists a sequence $v\geq v_{n}\downarrow0$
satisfying $\inf_{n\geq1}f_{\mathrm{s}}(v_{n})\geq f_{\mathrm{s}}(v)-\varepsilon$.
Denote $w_{n}=\max(u,v_{n})$. For any $\eta>0$, we can find $n_{0}\geq1$
such that $I_{\Phi}(v_{n})\leq\eta$ for every $n\geq n_{0}$. From
\[
I_{\Phi}(w_{n})\leq I_{\Phi}(u)+I_{\Phi}(v_{n})\leq1+\eta,
\]
it follows that $\Vert w_{n}\Vert_{\Phi}\leq1+\eta$, for every $n\geq n_{0}$.
Hence, for any $n\geq n_{0}$, we can write 
\begin{align*}
(1+\eta)\Vert f\Vert_{0} & \geq\Vert w_{n}\Vert_{\Phi}\Vert f\Vert_{0}\geq f(w_{n})=f_{\mathrm{c}}(w_{n})+f_{\mathrm{s}}(w_{n})\\
 & \geq f_{\mathrm{c}}(u)+f_{\mathrm{s}}(v_{n})\geq f_{\mathrm{c}}(u)+f_{\mathrm{s}}(v)-\varepsilon\\
 & \geq\Vert f_{\mathrm{c}}\Vert_{0}+\Vert f_{\mathrm{s}}\Vert_{0}-3\varepsilon.
\end{align*}
Since $\varepsilon,\eta>0$ are arbitrary, it follows that $\Vert f\Vert_{0}\geq\Vert f_{\mathrm{c}}\Vert_{0}+\Vert f_{\mathrm{s}}\Vert_{0}$. 
\end{proof}

\begin{prop}
\label{prop:functional_norm2} The norm of any functional $f=f_{v}+f_{\mathrm{s}}\in(L_{0}^{\Phi})^{*}$
can be expressed as 
\[
\Vert f\Vert=\inf\{\lambda>0:I_{\Phi^{*}}(v/\lambda)+\Vert f_{\mathrm{s}}/\lambda\Vert\leq1\}.
\]

\end{prop}

\begin{proof}
Without loss of generality, we assume that $\Vert f\Vert=1$ and $f\geq0$.
Take any $\lambda>0$ satisfying $I_{\Phi^{*}}(v/\lambda)+\Vert f_{\mathrm{s}}/\lambda\Vert\leq1$.
For any non-negative function $u\in L_{0}^{\Phi}$, and arbitrary
$k>0$ such that $I_{\Phi}(ku)<\infty$, we can write
\begin{align*}
\frac{1}{\lambda}|f(u)| & =\frac{1}{k}\Bigl(\int_{T}(ku)(v/\lambda)d\mu+f_{\mathrm{s}}(ku)/\lambda\Bigr)\\
 & \leq\frac{1}{k}(I_{\Phi}(ku)+I_{\Phi^{*}}(v/\lambda)+\Vert f_{\mathrm{s}}/\lambda\Vert)\\
 & \leq\frac{1}{k}(1+I_{\Phi}(ku))\leq\Vert u\Vert_{\Phi,0},
\end{align*}
where the first inequality follows from Young's inequality and expression
\prettyref{eq:norm_sing_func_2}. Thus we can conclude that
\begin{equation}
\Vert f\Vert\leq\inf\{\lambda>0:I_{\Phi^{*}}(v/\lambda)+\Vert f_{\mathrm{s}}/\lambda\Vert\leq1\}.\label{eq:functional_norm2_eq0}
\end{equation}
The function $v$ satisfies the inequality $I_{\Phi^{*}}(v)\leq1$.
This is a consequence of $\Vert v\Vert_{\Phi^{*}}=\Vert f_{v}\Vert\leq1$,
since $f_{v}(u)\leq f(u)\leq1$ for every non-negative function $u\in L_{0}^{\Phi}$
such that $\Vert u\Vert_{\Phi,0}\leq1$. According to \prettyref{lem:seq_un_IPhiS},
there exists a sequence of non-negative measurable functions $\{v_{n}\}$
such that $v_{n}\uparrow v$, and $I_{\Phi}((\Phi^{*})_{-}'(t,v_{n}(t)))<\infty$
for all $n\geq1$. Supposing that the first inequality in \prettyref{eq:functional_norm2_eq0}
is strict, we take some $\delta>0$ such that $I_{\Phi^{*}}(v)+\Vert f_{\mathrm{s}}\Vert>1+\delta$.
In view of \prettyref{eq:norm_sing_func_2}, we can find a non-negative
function $w\in\tilde{L}^{\Phi}$ such that $f_{\mathrm{s}}(w)\geq\Vert f_{\mathrm{s}}\Vert-\delta/2$.
Thus, from \prettyref{eq:f_s}, there exists a sequence $\{w_{n}\}$
satisfying $w\geq w_{n}\downarrow0$ and $\inf_{n\geq1}f_{\mathrm{s}}(w_{n})\geq f_{\mathrm{s}}(w)-\delta/4$.
Select a sufficiently large $n_{0}\geq1$ so that $I_{\Phi^{*}}(v_{n_{0}})\geq I_{\Phi^{*}}(v)-\delta/8$.
Denoting $u_{n_{0}}(t)=(\Phi^{*})_{-}'(t,v_{n_{0}}(t))$, we take
some $n_{1}\geq1$ for which the function $\widetilde{u}=\max(w_{n_{1}},u_{n_{0}})$
satisfies $I_{\Phi}(\widetilde{u})\leq I_{\Phi}(u_{n_{0}})+\delta/8$.
By these choices, we can write
\begin{align*}
f(\widetilde{u}) & =\int_{T}\widetilde{u}vd\mu+f_{\mathrm{s}}(\widetilde{u})\\
 & \geq\int_{T}u_{n_{0}}v_{n_{0}}d\mu+f_{\mathrm{s}}(w_{n_{1}})\\
 & \geq I_{\Phi}(u_{n_{0}})+I_{\Phi^{*}}(v_{n_{0}})+f_{\mathrm{s}}(w)-\delta/4\\
 & \geq I_{\Phi}(\widetilde{u})-\delta/8+I_{\Phi^{*}}(v)-\delta/8+\Vert f_{\mathrm{s}}\Vert-\delta/2-\delta/4\\
 & =I_{\Phi}(\widetilde{u})+I_{\Phi^{*}}(v)+\Vert f_{\mathrm{s}}\Vert-\delta\\
 & >1+I_{\Phi}(\widetilde{u})\\
 & \geq\Vert\widetilde{u}\Vert_{\Phi,0},
\end{align*}
which implies that $\Vert f\Vert>1$. Therefore, the first inequality
in \prettyref{eq:functional_norm2_eq0} cannot be strict.
\end{proof}

\section{Main results\label{sec:main_results}}

In this section, we provide a characterization of support functionals
and smooth points in $L_{0}^{\Phi}$ , and, as a result, we give necessary
and sufficient conditions for the smoothness of $L_{0}^{\Phi}$.

\subsection{Support functionals}

The characterization of support functionals at a function $u\in L_{0}^{\Phi}$
depends on whether the set $K(u)$ is empty or not.

\begin{prop}
\label{prop:f_support_Ku_neq_empty} Let $u\in L_{0}^{\Phi}\setminus\{0\}$
be such that $K(u)\neq\emptyset$. Then $f=f_{v}+f_{\mathrm{s}}\in(L_{0}^{\Phi})^{*}$
is a support functional at $u$ if and only if, for any $k\in K(u)$,
\begin{itemize}
\item [\normalfont{(i)}] $I_{\Phi^{*}}(v)+\Vert f_{\mathrm{s}}\Vert=1$,
\item [\normalfont{(ii)}] $\Vert f_{\mathrm{s}}\Vert=f_{\mathrm{s}}(ku)$,
and
\item [\normalfont{(iii)}] $\operatorname{sgn}v(t)=\operatorname{sgn}u(t)$
and $|v(t)|\in\partial\Phi(t,|ku(t)|)$ for $\mu$-a.e.\ $t\in T$.
\end{itemize}
\end{prop}
\begin{proof}
Suppose that (i)--(iii) are satisfied. By (i), we have that $\Vert f\Vert\leq1$.
For any $k\in K(u)$, we can write 
\begin{align*}
f(u) & =\frac{1}{k}\biggl(\int_{T}kuvd\mu+f_{\mathrm{s}}(ku)\biggr)\\
 & =\frac{1}{k}(I_{\Phi}(ku)+I_{\Phi^{*}}(v)+\Vert f_{\mathrm{s}}\Vert)\\
 & =\frac{1}{k}(1+I_{\Phi}(ku))=\Vert u\Vert_{\Phi,0},
\end{align*}
which implies that $\Vert f\Vert=1$. Therefore, $f$ is a support
functional at $u$. Conversely, let $f=f_{v}+f_{\mathrm{s}}\in(L_{0}^{\Phi})^{*}$
be a support functional at $u$. Using the expression $\frac{1}{k}(1+I_{\Phi}(ku))=\Vert u\Vert_{\Phi,0}=f_{v}(u)+f_{\mathrm{s}}(u)$,
for any $k\in K(u)$, we can write 
\begin{align*}
1 & =f_{v}(ku)-I_{\Phi}(ku)+f_{\mathrm{s}}(ku)\\
 & \leq I_{\Phi}(ku)+I_{\Phi^{*}}(v)-I_{\Phi}(ku)+f_{\mathrm{s}}(ku)\\
 & =I_{\Phi^{*}}(v)+f_{\mathrm{s}}(ku)\\
 & \leq I_{\Phi^{*}}(v)+\Vert f_{\mathrm{s}}\Vert\leq1.
\end{align*}
Then we obtain (i) and (ii), and $f_{v}(ku)=I_{\Phi^{*}}(v)+I_{\Phi}(ku)$,
from which (iii) follows.
\end{proof}

\begin{prop}
\label{prop:f_support_Ku_eq_empty} Let $u\in L_{0}^{\Phi}\setminus\{0\}$
be such that $K(u)=\emptyset$. Then $f=f_{v}+f_{\mathrm{s}}\in(L_{0}^{\Phi})^{*}$
is a support functional at $u$ if and only if
\begin{itemize}
\item [\normalfont{(i)}] $I_{\Phi^{*}}(v)+\Vert f_{\mathrm{s}}\Vert\leq1$,
and
\item [\normalfont{(ii)}] $v\chi_{\operatorname{supp}u}=\operatorname{sgn}u\cdot b_{\Phi^{*}}\chi_{\operatorname{supp}u}$.
\end{itemize}
\end{prop}
\begin{proof}
The assumption $K(u)=\emptyset$ implies that $I_{\Phi}(\lambda u)<\infty$
for all $\lambda>0$, and $\Vert u\Vert_{\Phi,0}=\int_{T}|u|b_{\Phi^{*}}d\mu$.
It is clear that if (i)--(ii) are satisfied then $f$ is a support
functional at $u$. Conversely, let $f=f_{v}+f_{\mathrm{s}}\in(L_{0}^{\Phi})^{*}$
be a support functional at $u$. Condition (i) follows from \prettyref{prop:functional_norm2}.
Clearly, $|v|\leq b_{\Phi^{*}}$. Suppose that the set $\{t\in\operatorname{supp}u:|v(t)|<b_{\Phi^{*}}(t)\}$
has non-zero measure. From $I_{\Phi}(\lambda u)<\infty$ for all $\lambda>0$,
we obtain that $f_{\mathrm{s}}(u)=0$. Then we can write
\begin{align*}
f(u) & =f_{v}(u)+f_{\mathrm{s}}(u)=f_{v}(u)\\
 & =\int_{T}uvd\mu\leq\int_{T}|u|\,|v|d\mu\\
 & <\int_{T}|u|b_{\Phi^{*}}d\mu=\Vert u\Vert_{\Phi,0},
\end{align*}
contradicting the assumption that $f$ is a support functional at
$u$. Therefore, $v\chi_{\operatorname{supp}u}=\operatorname{sgn}u\cdot b_{\Phi^{*}}\chi_{\operatorname{supp}u}$.
\end{proof}

\begin{cor}
\label{cor:Ku_more_than_one} Let $u\in L_{0}^{\Phi}\setminus\{0\}$
for which the set $K(u)\neq\emptyset$ is composed by more than one
element. Then there exists only one support functional at $u$, which
is given by $f_{v}$ with $v(t)=\operatorname{sgn}u(t)\cdot\Phi_{+}'(t,|k_{u}^{*}u(t)|)$. \end{cor}
\begin{proof}
From (ii) in \prettyref{prop:f_support_Ku_neq_empty}, we conclude
that every support functional at $u$ is order continuous. By the
definitions of $k_{u}^{*}$ and $k_{u}^{**}$, it is clear that $I_{\Phi^{*}}(\Phi_{+}'(t,|ku(t)|))=1$
for each $k\in[k_{u}^{*},k_{u}^{**})$. Consequently, $\Phi_{+}'(t,|k_{u}^{*}u(t)|)=\Phi_{-}'(t,|k_{u}^{**}u(t)|)$
for $\mu$-a.e.\ $t\in T$. Therefore, $v(t)=\operatorname{sgn}u(t)\cdot\Phi_{+}'(t,|k_{u}^{*}u(t)|)$
is the unique function satisfying $I_{\Phi^{*}}(v)=1$, and such that
$\operatorname{sgn}v(t)=\operatorname{sgn}u(t)$ and $|v(t)|\in\partial\Phi(t,|ku(t)|)$
for each $k\in K(u)$.
\end{proof}

\subsection{Smooth points}

To find necessary and sufficient conditions for a function in $L_{0}^{\Phi}$
to be a smooth point, we need some preliminary lemmas. The following
result is adapted from \cite[Lemma 6]{Grzaslewicz:1992} and \cite[Lemma 5]{Zbaszyniak:1994}.

\begin{lem}
\label{lem:disjoint_sequence} If the function $u\in\tilde{L}^{\Phi}$
satisfies $I_{\Phi}(\lambda u)=\infty$ for any $\lambda>1$, then
there exist non-increasing sequences of measurable sets $\{A_{n}\}$
and $\{B_{n}\}$, converging to the empty set, such that $A_{n}\cap B_{n}=\emptyset$
and $I_{\Phi}(\lambda u\chi_{A_{n}})=I_{\Phi}(\lambda u\chi_{B_{n}})=\infty$
for any $\lambda>1$ and $n\geq1$.\end{lem}
\begin{proof}
The proof is divided into three cases.

\textit{Case}\enskip{}1.\ Suppose that the measurable set $E=\{t\in T:|u(t)|=b_{\Phi}(t)\}$
has positive measure $\mu(E)>0$. Let $\{A_{n}\}$ and $\{B_{n}\}$
be non-increasing sequences of measurable sets, converging to the
empty set, such that $A_{n}\cap B_{n}=\emptyset$ and satisfying $0<\mu(E\cap A_{n})$
and $0<\mu(E\cap B_{n})$. Clearly, for each $n\geq1$, we have that
$I_{\Phi}(\lambda u\chi_{A_{n}})=I_{\Phi}(\lambda u\chi_{B_{n}})=\infty$
for any $\lambda>1$.

\textit{Case}\enskip{}2.\ Assume that $|u|<b_{\Phi}$, and for any
$\lambda>1$, the measurable set $F_{\lambda}=\{t\in T:|\lambda u(t)|>b_{\Phi}(t)\}$
has positive measure $\mu(F_{\lambda})>0$. Let $\{\lambda_{n}\}$
be a decreasing sequence in $(1,\infty)$ satisfying $\lambda_{n}\downarrow1$.
For each $n\geq1$, denote $F_{n}=F_{\lambda_{n}}$. Clearly, $0<\mu(F_{n})\downarrow0$.
For each $n\geq1$, take disjoint, measurable sets $G_{n}$ and $H_{n}$,
whose union is $G_{n}\cup H_{n}=F_{n}\setminus F_{n+1}$, and such
that $\mu(G_{n})>0$ and $\mu(H_{n})>0$ if $\mu(F_{n}\setminus F_{n+1})>0$,
or $\mu(G_{n})=\mu(H_{n})=0$ if $\mu(F_{n}\setminus F_{n+1})=0$.
For each $n\geq1$, define the disjoint sets $A_{n}=\bigcup_{k=n}^{\infty}G_{k}$
and $B_{n}=\bigcup_{k=n}^{\infty}H_{k}$. Clearly, we have that $\mu(A_{n})>0$
and $\mu(B_{n})>0$, for every $n\geq1$. Take any $\lambda>1$ and
$n\geq1$. For a sufficiently large $n_{0}\geq n$ such that $\lambda\geq\lambda_{n_{0}}$,
it follows that
\[
I_{\Phi}(\lambda u\chi_{A_{n}})=\sum_{k=n}^{\infty}I_{\Phi}(\lambda u\chi_{G_{k}})\geq\sum_{k=n_{0}}^{\infty}I_{\Phi}(\lambda_{k}u\chi_{G_{k}})=\infty.
\]
Similarly, we have that $I_{\Phi}(\lambda u\chi_{B_{n}})=\infty$
for any $\lambda>1$ and $n\geq1$.

\textit{Case}\enskip{}3.\ Now suppose that $|\overline{\lambda}u|<b_{\Phi}$
for some $\overline{\lambda}>1$. Let $\{\lambda_{n}\}$ be a decreasing
sequence in $(1,\overline{\lambda})$ such that $\lambda_{n}\downarrow1$.
Let $\{T_{n}\}$ be a non-decreasing sequence of measurable sets such
that $0<\mu(T_{n})<\infty$ and $\mu(T\setminus\bigcup_{n=1}^{\infty}T_{n})=0$.
Define the measurable sets $E_{n}^{m}=\{t\in T_{m}:\Phi(t,|\lambda_{n}u(t)|)\leq m\}$,
for all $n,m\geq1$. Clearly, $\chi_{E_{n}^{m}}\uparrow1$ as $m\rightarrow\infty$,
for each $n\geq1$. In view of $I_{\Phi}(\lambda_{1}u)=\infty$, we
can find $m_{1}\geq1$ such that $F_{1}=E_{1}^{m_{1}}$ satisfies
$2\leq I_{\Phi}(\lambda_{1}u\chi_{F_{1}})\leq m_{1}\mu(T_{m_{1}})<\infty$.
Obviously, $I_{\Phi}(\lambda_{n}u\chi_{T\setminus F_{1}})=\infty$
for all $n>1$. Similarly, we can find $m_{2}>m_{1}$ such that defining
$F_{2}=E_{2}^{m_{2}}\cap(T\setminus F_{1})$ we get $F_{1}\cap F_{2}=\emptyset$
and $2\leq I_{\Phi}(\lambda_{2}u\chi_{F_{2}})\leq m_{2}\mu(T_{m_{2}})<\infty$.
Thus $I_{\Phi}(\lambda_{n}u\chi_{T\setminus(F_{1}\cup F_{2})})=\infty$
for all $n>2$. Repeating these steps we obtain a sequence $\{F_{n}\}$
of pairwise disjoint sets such that $2\leq I_{\Phi}(\lambda_{n}u\chi_{F_{n}})<\infty$
for all $n\geq1$. Since the measure $\mu$ is non-atomic, there exist
disjoint, measurable sets $G_{n}$ and $H_{n}$, whose union is $F_{n}=G_{n}\cup H_{n}$,
such that
\[
I_{\Phi}(\lambda_{n}u\chi_{G_{n}})=I_{\Phi}(\lambda_{n}u\chi_{H_{n}})=\frac{1}{2}I_{\Phi}(\lambda_{n}u\chi_{F_{n}})\geq1.
\]
For each $n\geq1$, define the disjoint sets $A_{n}=\bigcup_{k=n}^{\infty}G_{k}$
and $B_{n}=\bigcup_{k=n}^{\infty}H_{k}$. Take any $\lambda>1$ and
$n\geq1$. For some $n_{0}\geq n$ such that $\lambda\geq\lambda_{n_{0}}$,
we can write
\[
I_{\Phi}(\lambda u\chi_{A_{n}})=\sum_{k=n}^{\infty}I_{\Phi}(\lambda u\chi_{G_{k}})\geq\sum_{k=n_{0}}^{\infty}I_{\Phi}(\lambda_{k}u\chi_{G_{k}})=\infty.
\]
Analogously, we obtain that $I_{\Phi}(\lambda u\chi_{B_{n}})=\infty$
for any $\lambda>1$ and $n\geq1$.
\end{proof}

\begin{lem}
\label{lem:disjoint_sets_singular} If the function $u\in\tilde{L}^{\Phi}$
satisfies $I_{\Phi}(\lambda u)=\infty$ for any $\lambda>1$, then
there exist two purely singular functionals $s_{1}\neq s_{2}$ in
$(L_{0}^{\Phi})^{*}$, with norms $\Vert s_{1}\Vert=\Vert s_{2}\Vert=1$,
and such that $s_{1}(u)=s_{2}(u)=1$.\end{lem}
\begin{proof}
According to \prettyref{lem:disjoint_sequence}, there exist non-increasing
sequences of measurable sets $\{A_{n}\}$ and $\{B_{n}\}$, converging
to the empty set, such that $A_{n}\cap B_{n}=\emptyset$ and $I_{\Phi}(\lambda u\chi_{A_{n}})=I_{\Phi}(\lambda u\chi_{B_{n}})=\infty$
for any $\lambda>1$, and all $n\geq1$. Without loss of generality,
we assume that $I_{\Phi}(u\chi_{A_{n}})\leq1$ and $I_{\Phi}(u\chi_{B_{n}})\leq1$
for all $n\geq1$. By this assumption, it follows that $\Vert u\chi_{A_{n}}\Vert_{\Phi}=\Vert u\chi_{B_{n}}\Vert_{\Phi}=1$
for all $n\geq1$. Denote the subspaces
\begin{align*}
\mathcal{E}_{1} & =\{w\in L^{\Phi}:\operatorname{supp}w\in T\setminus A_{n}\text{ for some }n\geq1\},\\
\mathcal{E}_{2} & =\{w\in L^{\Phi}:\operatorname{supp}w\in T\setminus B_{n}\text{ for some }n\geq1\}.
\end{align*}
It is clear that 
\[
\inf\{\Vert u-w\Vert_{\Phi}:w\in\mathcal{E}_{1}\}=\inf_{n\geq1}\Vert u\chi_{A_{n}}\Vert_{\Phi}=1.
\]
Hence the function $u$ does not belong to the closure of $\mathcal{E}_{1}$.
Similarly, $u$ is not in the closure of $\mathcal{E}_{2}$. By the
Hahn--Banach Theorem, we can find functionals $s_{1},s_{2}\in(L_{0}^{\Phi})^{*}$,
with norms $\Vert s_{1}\Vert=\Vert s_{2}\Vert=1$, and satisfying
$s_{1}(u)=s_{2}(u)=1$ and 
\begin{align*}
s_{1}(w) & =0,\qquad\text{for every }w\in\mathcal{E}_{1},\\
s_{2}(w) & =0,\qquad\text{for every }w\in\mathcal{E}_{2}.
\end{align*}
Since $B_{n}\subseteq T\setminus A_{n}$, we have that $s_{1}(u\chi_{B_{n}})=0$
and $s_{2}(u\chi_{B_{n}})=1$. Hence $s_{1}\neq s_{2}$. Clearly,
the positive and negative parts of $s_{1}$ vanish on $\mathcal{E}_{1}$.
For any non-negative $w\in L^{\Phi}$, it follows that
\begin{align*}
((s_{1})_{\pm})_{\mathrm{c}}(w) & =\inf\{\sup\nolimits _{n\geq1}(s_{1})_{\pm}(w_{n}):0\leq w_{n}\uparrow w\}\\
 & \leq\sup\nolimits _{n\geq1}(s_{1})_{\pm}(w\chi_{T\setminus A_{n}})=0.
\end{align*}
Therefore, $s_{1}$ is purely singular. Analogously, we have that
$s_{2}$ is purely singular.
\end{proof}

\begin{prop}
\label{prop:smooth_point_Ku_neq_empty} Let $u\in L_{0}^{\Phi}\setminus\{0\}$
be such that $K(u)\neq\emptyset$. Then $u$ is a smooth point if
and only if at least one of the following conditions is satisfied:
\begin{itemize}
\item [\normalfont{(i)}] $I_{\Phi^{*}}(\Phi_{-}'(t,|k_{u}^{*}u(t)|))=1$,
or 
\item [\normalfont{(ii)}] $I_{\Phi^{*}}(\Phi_{+}'(t,|k_{u}^{*}u(t)|))=1$
and $I_{\Phi}(\lambda u)<\infty$ for some $\lambda>k_{u}^{*}$.
\end{itemize}
\end{prop}
\begin{proof}
Assume that $u\in L_{0}^{\Phi}\setminus\{0\}$ is a smooth point.
If both conditions (i) and (ii) are not satisfied, then at least one
of the following expressions holds:
\begin{alignat}{2}
I_{\Phi^{*}}(\Phi_{-}'(t,|k_{u}^{*}u(t)|)) & <1 & \qquad\text{and}\qquad & I_{\Phi^{*}}(\Phi_{+}'(t,|k_{u}^{*}u(t)|))>1,\label{eq:prop_smooth point_A}\\
\intertext{or}I_{\Phi^{*}}(\Phi_{-}'(t,|k_{u}^{*}u(t)|)) & <1 & \qquad\text{and}\qquad & I_{\Phi}(\lambda u)=\infty\text{ for all }\lambda>k_{u}^{*}.\label{eq:prop_smooth point_B}
\end{alignat}
If \prettyref{eq:prop_smooth point_A} is satisfied, then we can find
a finite-valued, measurable function $v$ such that $\operatorname{sgn}v(t)=\operatorname{sgn}u(t)$
and $|v(t)|\in\partial\Phi(t,|k_{u}^{*}u(t)|)$ for $\mu$-a.e.\ $t\in T$,
and $I_{\Phi^{*}}(v)>1$. Denoting $\delta=I_{\Phi^{*}}(\Phi_{-}'(t,|k_{u}^{*}u(t)|))<1$,
we take any $\eta\in(0,1-\delta]$ such that $I_{\Phi^{*}}(v)\geq1+\eta$.
Then we can write
\[
I_{\Phi^{*}}(v)-I_{\Phi^{*}}(\Phi_{-}'(t,|k_{u}^{*}u(t)|))\geq1-\delta+\eta.
\]
Because the measure $\mu$ is non-atomic, there exists a measurable
set $E$ such that 
\begin{equation}
I_{\Phi^{*}}(v\chi_{E})-I_{\Phi^{*}}(\Phi_{-}'(t,|k_{u}^{*}u(t)|)\chi_{E}(t))=1-\delta+\eta.\label{eq:prop_smooth_point_delta}
\end{equation}
In view of $\eta\in(0,1-\delta]$, we can find disjoint, measurable
sets $A,B\subset E$ such that
\begin{align}
I_{\Phi^{*}}(v\chi_{A})-I_{\Phi^{*}}(\Phi_{-}'(t,|k_{u}^{*}u(t)|)\chi_{A}(t)) & =\eta,\label{eq:prop_smooth_point_eta_A}\\
I_{\Phi^{*}}(v\chi_{B})-I_{\Phi^{*}}(\Phi_{-}'(t,|k_{u}^{*}u(t)|)\chi_{B}(t)) & =\eta.\label{eq:prop_smooth_point_eta_B}
\end{align}
Clearly, the intersection of $A$ or $B$ with the set $\{t\in T:v(t)>\Phi_{-}'(t,|k_{u}^{*}u(t)|)\}$
has non-zero measure. Thus the following functions are different:
\begin{align*}
v_{1}(t) & :=v(t)\chi_{E\setminus A}(t)+\Phi_{-}'(t,|k_{u}^{*}u(t)|)\chi_{T\setminus(E\setminus A)}(t),\\
v_{2}(t) & :=v(t)\chi_{E\setminus B}(t)+\Phi_{-}'(t,|k_{u}^{*}u(t)|)\chi_{T\setminus(E\setminus B)}(t).
\end{align*}
From \prettyref{eq:prop_smooth_point_delta} and \prettyref{eq:prop_smooth_point_eta_A},
we can write
\begin{align*}
I_{\Phi^{*}}(v_{1}) & =I_{\Phi^{*}}(v\chi_{E\setminus A})+I_{\Phi^{*}}(\Phi_{-}'(t,|k_{u}^{*}u(t)|)\chi_{T\setminus(E\setminus A)}(t))\\
 & =I_{\Phi^{*}}(v\chi_{E})-I_{\Phi^{*}}(v\chi_{A})+I_{\Phi^{*}}(\Phi_{-}'(t,|k_{u}^{*}u(t)|))-I_{\Phi^{*}}(\Phi_{-}'(t,|k_{u}^{*}u(t)|)\chi_{E\setminus A}(t))\\
 & =I_{\Phi^{*}}(v\chi_{E})-I_{\Phi^{*}}(\Phi_{-}'(t,|k_{u}^{*}u(t)|)\chi_{E}(t))+\delta-[I_{\Phi^{*}}(v\chi_{A})-I_{\Phi^{*}}(\Phi_{-}'(t,|k_{u}^{*}u(t)|)\chi_{A}(t))]\\
 & =1-\delta+\eta+\delta-\eta=1.
\end{align*}
Analogously, from \prettyref{eq:prop_smooth_point_delta} and \prettyref{eq:prop_smooth_point_eta_B},
it follows that $I_{\Phi_{2}}(v_{2})=1$. According to \prettyref{prop:f_support_Ku_neq_empty},
$f_{v_{1}}$ and $f_{v_{2}}$ are different support functionals at
$u$. 

Now suppose that \prettyref{eq:prop_smooth point_B} is satisfied.
Let $v(t)=\operatorname{sgn}u(t)\cdot\Phi_{-}'(t,|k_{u}^{*}u(t)|)$.
Using \prettyref{lem:disjoint_sets_singular}, we can find two purely
singular functionals $s_{1}\neq s_{2}$ in $(L_{0}^{\Phi})^{*}$,
with norms $\Vert s_{1}\Vert=\Vert s_{2}\Vert=1-I_{\Phi^{*}}(v)$,
and such that $s_{1}(u)=s_{2}(u)=1-I_{\Phi^{*}}(v)$. According to
the proof of \prettyref{cor:Ku_more_than_one}, if the set $K(u)$
is composed by more than one element, then $I_{\Phi^{*}}(v)=1$, which
is a contradiction to \prettyref{eq:prop_smooth point_B}. Consequently,
$K(u)=\{k_{u}^{*}\}$. Then we conclude that $f_{1}$ and $f_{2}$
satisfy conditions (i) and (ii) in \prettyref{prop:f_support_Ku_neq_empty},
which shows that $f_{1}$ and $f_{2}$ are different support functionals
at $u$. Therefore, if $u\in L^{\Phi}\setminus\{0\}$ is a smooth
point, then at least one of conditions (i) or (ii) holds.

Next we will show that if at least one of conditions (i) or (ii) is
satisfied, then $u$ is a smooth point. Let $v$ by any measurable
function such that $\operatorname{sgn}v(t)=\operatorname{sgn}u(t)$
and $|v(t)|\in\partial\Phi(t,|k_{u}^{*}u(t)|)$ for $\mu$-a.e.\ $t\in T$.
Assume that (i) is satisfied. It is clear that $I_{\Phi^{*}}(v)\geq1$,
and that $I_{\Phi^{*}}(v)=1$ if and only if $v(t)=v_{0}(t)=\operatorname{sgn}u(t)\cdot\Phi_{-}'(t,|k_{u}^{*}u(t)|)$.
In view of \prettyref{prop:f_support_Ku_neq_empty}, $f_{v_{0}}$
is the unique support functional at $u$. Now suppose that (ii) holds.
Then does not exist a non-zero, purely singular functional $s$ such
that $\Vert s\Vert=s(k_{u}^{*}u)$, since in this case we can write
that $0<\Vert s\Vert=s(k_{u}^{*}u)<s(\lambda u)\leq\Vert s\Vert$
for some $\lambda>k_{u}^{*}$ such that $I_{\Phi}(\lambda u)<\infty$.
From (ii) in \prettyref{prop:f_support_Ku_neq_empty}, it follows
that every support functional at $u$ is order continuous. Assume
that $f_{v}$ is a support functional at $u$. Hence $I_{\Phi^{*}}(v)=1$.
If $k_{u}^{*}=k_{u}^{**}$ then $I_{\Phi^{*}}(v)=1$ implies that
$v(t)=v_{1}(t)=\operatorname{sgn}u(t)\cdot\Phi_{+}'(t,|k_{u}^{*}u(t)|)$.
In view of \prettyref{cor:Ku_more_than_one}, if the set $K(u)$ is
composed by more than one element, then $u$ is a smooth point and
$v(t)=v_{1}(t)=\operatorname{sgn}u(t)\cdot\Phi_{+}'(t,|k_{u}^{*}u(t)|)$.
Thus, assuming that (ii) is satisfied, we have that $f_{v_{1}}$ is
the unique support functional at $u$. 
\end{proof}

\begin{prop}
\label{prop:smooth_point_Ku_eq_empty} Let $u\in L_{0}^{\Phi}\setminus\{0\}$
be such that $K(u)=\emptyset$. Then $u$ is a smooth point if and
only if
\begin{itemize}
\item [\normalfont{(i)}] $I_{\Phi^{*}}(b_{\Phi^{*}}\chi_{\operatorname{supp}u})=1$
and $a_{\Phi^{*}}\chi_{T\setminus\operatorname{supp}u}=0$, or
\item [\normalfont{(ii)}] $I_{\Phi^{*}}(b_{\Phi^{*}})<1$ and $\mu(T\setminus\operatorname{supp}u)=0$.
\end{itemize}
\end{prop}
\begin{proof}
Assume that $u$ is a smooth point. Let $f=f_{v}+f_{\mathrm{s}}\in(L_{0}^{\Phi})^{*}$
be a support functional at $u$. According to \prettyref{prop:f_support_Ku_eq_empty},
we have that $I_{\Phi^{*}}(v)+\Vert f_{\mathrm{s}}\Vert\leq1$ and
$v\chi_{\operatorname{supp}u}=\operatorname{sgn}u\cdot b_{\Phi^{*}}\chi_{\operatorname{supp}u}$.
Suppose that $\Phi$ does not satisfy the $\Delta_{2}$-condition.
If $I_{\Phi^{*}}(v)<1$ then any functional $g=f_{v}+s$ such that
$I_{\Phi^{*}}(v)+\Vert s\Vert\leq1$, where $s$ is a purely singular
functional, is a support functional at $u$. Consequently, $I_{\Phi^{*}}(v)=1$.
We cannot have $I_{\Phi^{*}}(b_{\Phi^{*}}\chi_{\operatorname{supp}u})<1$,
since $f_{w}$ with $w=\operatorname{sgn}u\cdot b_{\Phi^{*}}\chi_{\operatorname{supp}u}\neq v$
would be a support functional at $u$. By $I_{\Phi^{*}}(v\chi_{T\setminus\operatorname{supp}u})=0$,
it follows that $0\leq|v|\chi_{T\setminus\operatorname{supp}u}\leq a_{\Phi^{*}}\chi_{T\setminus\operatorname{supp}u}$.
If $a_{\Phi^{*}}\chi_{T\setminus\operatorname{supp}u}>0$ then it
is clear that $f_{w}$ with $w=v\chi_{\operatorname{supp}u}+\frac{1}{2}(a_{\Phi^{*}}-|v|)\chi_{T\setminus\operatorname{supp}u}\neq v$
is a support functional at $u$. Hence $a_{\Phi^{*}}\chi_{T\setminus\operatorname{supp}u}=0$.
Therefore, (i) holds. Now suppose that $\Phi$ satisfies the $\Delta_{2}$-condition.
Thus $f_{\mathrm{s}}=0$. If $I_{\Phi^{*}}(b_{\Phi^{*}}\chi_{\operatorname{supp}u})=1$
then proceeding as above we obtain that $a_{\Phi^{*}}\chi_{T\setminus\operatorname{supp}u}=0$.
Assume that $I_{\Phi^{*}}(b_{\Phi^{*}}\chi_{\operatorname{supp}u})<1$.
If $\mu(T\setminus\operatorname{supp}u)\neq0$ then it is clear that
we can find $w\in\tilde{L}^{\Phi}$ with $w\neq v$ and $I_{\Phi^{*}}(w)\leq1$.
Hence $f_{w}$ is a support functional at $u$. Consequently, $\mu(T\setminus\operatorname{supp}u)=0$.
Then (i) or (ii) is satisfied.

Conversely, if (i) holds then in view of \prettyref{prop:f_support_Ku_eq_empty}
it is clear that $f_{v}$ with $v=b_{\Phi^{*}}\chi_{\operatorname{supp}u}$
is the unique support functional at $u$. Assume that (ii) is satisfied.
For any $w\in L^{\Phi}$, and $\lambda>0$, we have that 
\begin{align*}
I_{\Phi}(\lambda w) & \leq I_{\Phi}(\lambda w)+I_{\Phi^{*}}(\Phi_{-}'(t,|\lambda w(t)|))\\
 & =\int_{T}|\lambda w(t)|\Phi_{-}'(t,|\lambda w(t)|)d\mu\\
 & \leq\int_{T}|\lambda w|b_{\Phi^{*}}d\mu\\
 & =\Vert\lambda w\Vert_{\Phi,0}<\infty.
\end{align*}
Then $\Phi$ satisfies the $\Delta_{2}$-condition, and any functional
in $(L^{\Phi})^{*}$ is order continuous. Consequently, $f_{v}$ with
$v=b_{\Phi^{*}}$ is the unique support functional at $u$.
\end{proof}

\subsection{Smoothness of $L_{0}^{\Phi}$}

Below we state the main result of this paper, which provides necessary
and sufficient conditions for the smoothness of $L_{0}^{\Phi}$.

\begin{prop}
\label{prop:smoothness} The Musielak--Orlicz space $L_{0}^{\Phi}$
is smooth if and only if
\begin{itemize}
\item [\normalfont{(a)}] $\Phi^{*}(t,b_{\Phi^{*}}(t))=\infty$ for $\mu$-a.e.\ $t\in T$,
\item [\normalfont{(b)}] $\Phi$ satisfies the $\Delta_{2}$-condition,
and
\item [\normalfont{(c)}] $\Phi(t,\cdot)$ is continuously differentiable
(with $\Phi_{+}'(t,0)=0$) for $\mu$-a.e.\ $t\in T$.
\end{itemize}
\end{prop}

The proof of this proposition requires some preliminary results and
observations. Let $u\in L_{0}^{\Phi}\setminus\{0\}$ be such that
$K(u)\neq\emptyset$. Suppose that $f=f_{v}+f_{\mathrm{s}}$ is a
support functional of $u$ with non-zero singular component $f_{\mathrm{s}}\neq0$.
From condition (ii) in \prettyref{prop:f_support_Ku_neq_empty}, we
have that $\Vert f_{\mathrm{s}}\Vert=\sup_{u\in\tilde{L}^{\Phi}}|f(u)|=f_{\mathrm{s}}(ku)$
for any $k\in K(u)$. This result implies that $K(u)=\{k\}$, and
$I_{\Phi}(\lambda u)=\infty$ for any $\lambda>k$. Moreover, in view
of (ii) in \prettyref{prop:f_support_Ku_neq_empty}, it follows that
$I_{\Phi^{*}}(v)=1-\Vert f_{\mathrm{s}}\Vert<1$. Thus the existence
of a function $u\in L_{0}^{\Phi}$ such that $I_{\Phi}(\lambda u)=\infty$
for any $\lambda>k$, and $I_{\Phi^{*}}(\Phi_{-}'(t,|ku(t)|))<1$,
is a necessary condition for the existence of a support functional
with non-zero singular component. Thanks to the result below, we can
prove \prettyref{prop:smoothness} without using the Bishop--Phelps
Theorem (cf.~\cite[Theorem~2.3]{Hudzik:1997}).

\begin{prop}
\label{prop:IPhiStarLessOne} Let $\Phi$ be a Musielak--Orlicz function
not satisfying the $\Delta_{2}$-condition. Then there exists a function
$u\in\tilde{L}^{\Phi}$ such that $I_{\Phi}(\lambda u)=\infty$ for
any $\lambda>1$, and $I_{\Phi^{*}}(\Phi_{+}'(t,|u(t)|))<1$.
\end{prop}

Notice that, for the function $u\in\tilde{L}^{\Phi}$ in the proposition
above, we have that $K(u)=\{1\}$. To show \prettyref{prop:IPhiStarLessOne},
we make use of the lemma below, which is stated without proof (see
\cite[Lemma 8.3]{Musielak:1983}).

\begin{lem}
\label{lem:seq} Let $\mu$ be a non-atomic, $\sigma$-finite measure.
If $\{\alpha_{n}\}$ is a sequence of positive, real numbers, and
$\{u_{n}\}$ is a sequence of finite-valued, non-negative, measurable
functions, such that 
\[
\int_{T}u_{n}d\mu\geq2^{n}\alpha_{n},\quad\text{for all }n\geq1,
\]
then there exist an increasing sequence $\{n_{i}\}$ of natural numbers
and a sequence $\{A_{i}\}$ of pairwise disjoint, measurable sets
such that 
\[
\int_{A_{i}}u_{n_{i}}d\mu=\alpha_{n_{i}},\quad\text{for all }i\geq1.
\]

\end{lem}

One can easily verify that the $\Delta_{2}$-condition given by \prettyref{eq:D2_0}
is equivalent to the existence of a constant $\alpha>0$, and a non-negative
function $f\in\tilde{L}^{\Phi}$ such that
\begin{equation}
\alpha\Phi(t,u)\leq\Phi(t,\tfrac{1}{2}u),\quad\text{for all }u>f(t).\label{eq:D2_1}
\end{equation}
Moreover, a Musielak--Orlicz function $\Phi$ satisfies the $\Delta_{2}$-condition
if, and only if, for every $\lambda\in(0,1)$, there exist a constant
$\alpha_{\lambda}\in(0,1)$, and a non-negative function $f_{\lambda}\in\tilde{L}^{\Phi}$
such that 
\begin{equation}
\alpha_{\lambda}\Phi(t,u)\leq\Phi(t,\lambda u),\quad\text{for all }u>f_{\lambda}(t).\label{eq:D2_2}
\end{equation}
We will use this observation to prove the next result.

\begin{lem}
\label{lem:technical} Let $\Phi$ be a Musielak--Orlicz function
not satisfying the $\Delta_{2}$-condition. Assume that $\Phi(t,b_{\Phi}(t))=\infty$
for $\mu$-a.e.\ $t\in T$. Then we can find a strictly increasing
sequence $\{\lambda_{n}\}$ in $(0,1)$ converging upward to $1$,
and sequences $\{u_{n}\}$ and $\{A_{n}\}$ of finite-valued, measurable
functions, and pairwise disjoint, measurable sets, respectively, such
that
\begin{equation}
I_{\Phi}(u_{n}\chi_{A_{n}})=1\quad\text{and}\quad I_{\Phi}(\lambda_{n}u_{n}\chi_{A_{n}})\leq2^{-n},\quad\text{for all }n\geq1.\label{eq:u_n_lambda_n}
\end{equation}
\end{lem}
\begin{proof}
Because the Musielak--Orlicz function $\Phi$ does not satisfy the
$\Delta_{2}$-condition, for any $\lambda\in(0,1)$, there do not
exist a constant $\alpha\in(0,1)$ and a non-negative function $f\in\tilde{L}^{\Phi}$
such that 
\begin{equation}
\alpha\Phi(t,u)\leq\Phi(t,\lambda u),\quad\text{for all }u>f(t).\label{eq:class_ineq2}
\end{equation}
Let $\{\lambda'_{m}\}$ be a strictly increasing sequence in $(0,1)$
such that $\lambda'_{m}\uparrow1$. Define the non-negative functions
\[
f_{m}(t)=\sup\{u\in(0,b_{\Phi}(t)):2^{-m}\Phi(t,u)>\Phi(t,\lambda'_{m}u)\},\quad\text{for all }m\geq1,
\]
where we adopt the convention that $\sup\emptyset=0$. Since \prettyref{eq:class_ineq2}
is not satisfied, we have that $I_{\Phi}(f_{m})=\infty$ for each
$m\geq1$. For every rational number $r>0$, define the measurable
sets 
\[
A_{m,r}=\{t\in T:r\in(0,b_{\Phi}(t))\text{ and }2^{-m}\Phi(t,u)>\Phi(t,\lambda'_{m}u)\}
\]
and the simple functions $u_{m,r}=r\chi_{A_{m,r}}$. For $r=0$, set
$u_{m,r}=0$. Let $\{r_{i}\}$ be an enumeration of the non-negative
rational numbers with $r_{1}=0$. Define the non-negative, simple
functions $v_{m,k}=\max_{1\leq i\leq k}u_{m,r_{i}}$, for each $m,k\geq1$.
By the left-continuity of $\Phi(t,\cdot)$, it follows that $v_{m,k}\uparrow f_{m}$
as $k\rightarrow\infty$. In virtue of the Monotone Convergence Theorem,
for each $m\geq1$, we can find some $k_{m}\geq1$ such that the function
$v_{m}=v_{m,k_{m}}$ satisfies $I_{\Phi}(v_{m})\geq2^{m}$. Clearly,
we have that $\Phi(t,v_{m}(t))<\infty$ and $2^{-m}\Phi(t,v_{m}(t))\geq\Phi(t,\lambda'_{m}v_{m}(t))$.
By \prettyref{lem:seq}, there exist an increasing sequence $\{m_{n}\}$
of indices and a sequence $\{A_{n}\}$ of pairwise disjoint, measurable
sets such that $I_{\Phi}(v_{m_{n}}\chi_{A_{n}})=1$. Taking $\lambda_{n}=\lambda'_{m_{n}}$,
$u_{n}=v_{m_{n}}$ and $A_{n}$, we obtain \prettyref{eq:u_n_lambda_n}.
\end{proof}

\begin{proof}[Proof of \prettyref{prop:IPhiStarLessOne}]
 Suppose that the measurable set $E=\{t\in T:\Phi(t,b_{\Phi}(t))<\infty\}$
has positive measure $\mu(E)>0$. Take a measurable set $F\subseteq E$
such that $\mu(F)>0$ and $I_{\Phi}(b_{\Phi}\chi_{F})<\infty$. Since
the measure $\mu$ is non-atomic, we can find pairwise disjoint, measurable
sets $A_{n}\subset F$ such that $\mu(A_{n})>0$ and $F=\bigcup_{n=1}^{\infty}A_{n}$.
Let $\{\lambda_{n}\}$ be a strictly increasing sequence in $(0,1)$
such that $\lambda_{n}\uparrow1$. For each $n\geq1$, select a measurable
set $B_{n}\subseteq A_{n}$ such that $\mu(B_{n})>0$ and 
\begin{equation}
I_{\Phi^{*}}(\Phi_{+}'(t,\lambda_{n}b_{\Phi}(t)\chi_{B_{n}}(t)))<2^{-n}.\label{eq:propIPhiStarLessOne}
\end{equation}
Define the function $u=\sum_{n=1}^{\infty}\lambda_{n}b_{\Phi}\chi_{B_{n}}$.
Clearly, $I_{\Phi}(u)<\infty$. For any $\lambda>1$, and some $n_{0}\geq1$
such that $\lambda\lambda_{n_{0}}>1$, we have that
\[
I_{\Phi}(\lambda u)=\sum_{n=1}^{\infty}I_{\Phi}(\lambda\lambda_{n}b_{\Phi}\chi_{B_{n}})\geq I_{\Phi}(\lambda\lambda_{n_{0}}b_{\Phi}\chi_{B_{n_{0}}})=\infty.
\]
By \prettyref{eq:propIPhiStarLessOne}, it follows that
\[
I_{\Phi^{*}}(\Phi_{+}'(t,u(t)))=\sum_{n=1}^{\infty}I_{\Phi^{*}}(\Phi_{+}'(t,\lambda_{n}b_{\Phi}(t)\chi_{B_{n}}(t)))<\sum_{n=1}^{\infty}2^{-n}=1.
\]

Now assume that $\Phi(t,b_{\Phi}(t))=\infty$ for $\mu$-a.e.\ $t\in T$.
Let $\{\lambda_{n}\}$, $\{u_{n}\}$ and $\{A_{n}\}$ be the sequences
in the statement of \prettyref{lem:technical}. For a sufficiently
large natural number $n_{0}>1$ such that $\lambda_{n_{0}}>1/n_{0}$
and $\sum_{n=n_{0}}^{\infty}n2^{-n}<1$, we define the function $u=\sum_{n=n_{0}}^{\infty}\lambda_{n}'u_{n}\chi_{A_{n}}$,
where $\lambda_{n}'=\lambda_{n}-1/n$, for each $n\geq n_{0}$. Then
we can write
\[
I_{\Phi}(u)=\sum_{n=n_{0}}^{\infty}I_{\Phi}(\lambda_{n}'u_{n}\chi_{A_{n}})\leq\sum_{n=n_{0}}^{\infty}I_{\Phi}(\lambda_{n}u_{n}\chi_{A_{n}})\leq\sum_{n=n_{0}}^{\infty}2^{-n}<\infty.
\]
Given any $\lambda>1$, we take some $n_{1}\geq n_{0}$ satisfying
$\lambda\lambda_{n_{1}}'\geq1$, so that
\[
I_{\Phi}(\lambda u)=\sum_{n=n_{0}}^{\infty}I_{\Phi}(\lambda\lambda_{n}'u_{n}\chi_{A_{n}})\geq\sum_{n=n_{1}}^{\infty}I_{\Phi}(u_{n}\chi_{A_{n}})=\infty.
\]
For each $n\geq n_{0}$, we obtain that 
\begin{align*}
I_{\Phi^{*}}(\Phi_{+}'(t,\lambda_{n}'u_{n}(t)\chi_{A_{n}}(t))) & \leq I_{\Phi}(\lambda_{n}'u_{n}\chi_{A_{n}})+I_{\Phi^{*}}(\Phi_{+}'(t,\lambda_{n}'u_{n}(t)\chi_{A_{n}}(t)))\\
 & =\int_{A_{n}}\lambda_{n}'u_{n}(t)\Phi_{+}'(t,\lambda_{n}'u_{n}(t))d\mu\\
 & \leq\lambda_{n}'\frac{1}{\lambda_{n}-\lambda_{n}'}[I_{\Phi}(\lambda_{n}u_{n}\chi_{A_{n}})-I_{\Phi}(\lambda_{n}'u_{n}\chi_{A_{n}})]\\
 & \leq\frac{1}{\lambda_{n}-\lambda_{n}'}I_{\Phi}(\lambda_{n}u_{n}\chi_{A_{n}})\\
 & \leq n2^{-n}.
\end{align*}
Hence it follows that
\[
I_{\Phi^{*}}(\Phi_{+}'(t,u(t)))=\sum_{n=n_{0}}^{\infty}I_{\Phi^{*}}(\Phi_{+}'(t,\lambda_{n}'u_{n}(t)\chi_{A_{n}}(t)))\leq\sum_{n=n_{0}}^{\infty}n2^{-n}<1,
\]
which completes the proof.
\end{proof}

\begin{rem}
\label{rem:u_stars} Let $\Phi$ be a Musielak--Orlicz function not
satisfying the $\Delta_{2}$-condition and such that $\Phi(t,b_{\Phi}(t))=\infty$
for $\mu$-a.e.\ $t\in T$. Then we can find functions $u_{*}$ and
$u^{*}$ in $L^{\Phi}$ such that 
\begin{equation}
\left\{ \begin{array}{ll}
I_{\Phi}(\lambda u_{*})<\infty, & \quad\text{for }0\leq\lambda\leq1,\\
I_{\Phi}(\lambda u_{*})=\infty, & \quad\text{for }1<\lambda,
\end{array}\right.\label{eq:substar}
\end{equation}
and
\begin{equation}
\left\{ \begin{array}{ll}
I_{\Phi}(\lambda u^{*})<\infty, & \quad\text{for }0\leq\lambda<1,\\
I_{\Phi}(\lambda u^{*})=\infty, & \quad\text{for }1\leq\lambda.
\end{array}\right.\label{eq:upstar}
\end{equation}
We construct these functions using the sequences $\{\lambda_{n}\}$,
$\{u_{n}\}$ and $\{A_{n}\}$ in \prettyref{lem:technical}. Define
$u_{*}=\sum_{n=1}^{\infty}\lambda_{n}u_{n}\chi_{A_{n}}$ and $u^{*}=\sum_{n=1}^{\infty}u_{n}\chi_{A_{n}}$.
For any $0\leq\lambda\leq1$, we have that 
\[
I_{\Phi}(\lambda u_{*})\leq I_{\Phi}(u_{*})=\sum_{n=1}^{\infty}I_{\Phi}(\lambda_{n}u_{n}\chi_{A_{n}})\leq\sum_{n=1}^{\infty}2^{-n}=1.
\]
For any $\lambda>1$, take a natural number $n_{0}\geq1$ such that
$\lambda\lambda_{n_{0}}\geq1$. Then we can write 
\[
I_{\Phi}(\lambda u_{*})=\sum_{n=1}^{\infty}I_{\Phi}(\lambda\lambda_{n}u_{n}\chi_{A_{n}})\geq\sum_{n=n_{0}}^{\infty}I_{\Phi}(u_{n}\chi_{A_{n}})=\infty.
\]
With respect to $u^{*}$, it is clear that $I_{\Phi}(\lambda u^{*})=\infty$
for any $\lambda\geq1$. If $\lambda<1$ and the natural number $n_{0}\geq1$
is such that $\lambda\leq\lambda_{n_{0}}$, we obtain that
\[
I_{\Phi}(\lambda u^{*})=\sum_{n=1}^{\infty}I_{\Phi}(\lambda u_{n}\chi_{A_{n}})\leq\sum_{n=1}^{n_{0}-1}I_{\Phi}(\lambda u_{n}\chi_{A_{n}})+\sum_{n=n_{0}}^{\infty}I_{\Phi}(\lambda_{n}u_{n}\chi_{A_{n}})<\infty.
\]
Thus the functions $u_{*}$ and $u^{*}$ satisfy \prettyref{eq:substar}
and \prettyref{eq:upstar}.
\end{rem}

Thanks to the lemma below, we avoid the use of measurable selectors
in the proof of \prettyref{prop:smoothness} (cf.~\cite[Theorem~2.2]{Hudzik:1997}).

\begin{lem}
\label{lem:u_delta} Let $\Phi$ be a finite-valued Musielak--Orlicz
function. For any $\delta>0$, the function 
\[
u_{\delta}(t)=\sup\{u\geq0:\Phi_{+}'(t,x)-\Phi_{-}'(t,x)<\delta\text{ for all }0\leq x\leq u\}
\]
is measurable (where we adopted the convention $\Phi_{-}'(t,0)=0$).
Moreover, denoting $H=\{t\in T:u_{\delta}(t)<\infty\}$, then $\Phi_{+}'(t,u_{\delta}(t))-\Phi_{-}'(t,u_{\delta}(t))\geq\delta$
for $\mu$-a.e.\ $t\in H$.\end{lem}
\begin{proof}
Without loss of generality, we assume that, for each $t\in T$, the
function $\Phi(t,\cdot)$ satisfies conditions (i) and (ii) in the
definition of Musielak--Orlicz functions. For any $\varepsilon>0$,
we define the function 
\[
u_{\delta,\varepsilon}(t)=\sup\{u\geq0:\Phi_{-}'(t,x+\varepsilon)-\Phi_{-}'(t,x)<\delta\text{ for all }0\leq x\leq u\}.
\]
We will verify that $u_{\delta,\varepsilon}$ is measurable. Fixed
any $u\geq0$, denote $A_{u}=\{t\in T:\Phi_{-}'(t,x+\varepsilon)-\Phi_{-}'(t,x)<\delta\text{ for all }0\leq x\leq u\}$.
Let $\{\delta_{n}\}$ be a sequence in $(0,\delta)$ such that $\delta_{n}\uparrow\delta$.
Now define $\widetilde{A}_{u}=\bigcup_{n=1}^{\infty}\bigcap_{r\in[0,u]\cap\mathbb{Q}}B_{r}^{n}$,
where $B_{r}^{n}=\{t\in T:\Phi_{-}'(t,r+\varepsilon)-\Phi_{-}'(t,r)<\delta_{n}\}$,
for each $r\in[0,u]\cap\mathbb{Q}$ and $n\geq1$. Clearly, the sets
$B_{r}^{n}$ are measurable, and then $\widetilde{A}_{u}$ is also
measurable. We will show that $A_{u}$ and $\widetilde{A}_{u}$ coincide.
It is clear that $A_{u}\subseteq\widetilde{A}_{u}$. Fix any $t\in\widetilde{A}_{u}$.
Hence $t\in\bigcap_{r\in[0,u]\cap\mathbb{Q}}B_{r}^{n_{0}}$ for some
$n_{0}\geq1$. For any $x\in(0,u]$, take a sequence of rational numbers
$\{r_{k}\}$ in $(0,x)$ such that $r_{k}\uparrow x$. Since $\Phi_{-}'(t,r_{k}+\varepsilon)-\Phi_{-}'(t,r_{k})<\delta_{n_{0}}<\delta$
for all $k\geq1$, it follows that $\Phi_{-}'(t,x+\varepsilon)-\Phi_{-}'(t,x)\leq\delta_{n_{0}}<\delta$.
Thus $t\in A_{u}$, which implies that $A_{u}=\widetilde{A}_{u}$.
Consequently, $A_{u}$ is measurable. Let $\{r_{i}\}$ be an enumeration
of the non-negative, rational numbers with $r_{1}=0$. Clearly, the
non-negative simple functions $\overline{u}_{n}=\max_{1\leq i\leq n}r_{i}\chi_{A_{r_{i}}}$
converge upward to $u_{\delta,\varepsilon}$. Thus the function $u_{\delta,\varepsilon}$
is measurable.

If we show that $u_{\delta,\varepsilon}\uparrow u_{\delta}$ as $\varepsilon\downarrow0$,
then $u_{\delta}$ is measurable. It is clear that $u_{\delta,\varepsilon_{2}}(t)\leq u_{\delta,\varepsilon_{1}}(t)\leq u_{\delta}(t)$
for any $0<\varepsilon_{1}\leq\varepsilon_{2}$. Hence $u_{\delta,\varepsilon}(t)\uparrow c(t)$
for some measurable function $c$. Let $t\in T$ be such that $c(t)=0$.
In view of $u_{\delta,\varepsilon}(t)\leq c(t)=0$ for all $\varepsilon>0$,
we obtain that $\Phi_{-}'(t,0+\varepsilon)=\Phi_{-}'(t,0+\varepsilon)-\Phi_{-}'(t,0)\geq\delta$
for all $\varepsilon>0$. Thus $\Phi_{+}'(t,0)-\Phi_{-}'(t,0)\geq\delta$,
which implies that $u_{\delta}(t)=0$. Now suppose that there exists
some $t\in T$ satisfying $0<c(t)<u_{\delta}(t)$. Clearly, $\Phi_{+}'(t,c(t))-\Phi_{-}'(t,c(t))<\delta$.
Then we can find $\eta>0$ so that
\begin{equation}
\Phi_{-}'(t,x_{2})-\Phi_{-}'(t,x_{1})<\delta,\qquad\text{for all }x_{1}\in(c(t)-\eta,c(t)],x_{2}\in(c(t),c(t)+\eta).\label{eq:u_delta}
\end{equation}
Let $\varepsilon\in(0,\eta/2)$ be such that $u_{\delta,\varepsilon}(t)\in(c(t)-\eta/2,c(t)]$.
By the definition of $u_{\delta,\varepsilon}(t)$, there exists $x_{0}\in[u_{\delta,\varepsilon}(t),u_{\delta,\varepsilon}(t)+\eta/2)$
with $\Phi_{-}'(t,x_{0}+\varepsilon)-\Phi_{-}'(t,x_{0})\geq\delta$.
Denoting $x_{1}=\min(x_{0},c(t))$ and $x_{2}=x_{0}+\varepsilon$,
we can write $\Phi_{-}'(t,x_{2})-\Phi_{-}'(t,x_{1})\geq\Phi_{-}'(t,x_{0}+\varepsilon)-\Phi_{-}'(t,x_{0})\geq\delta$.
This is a contradiction to \prettyref{eq:u_delta}, since $x_{1}\in(c(t)-\eta,c(t)]$
and $x_{2}\in(c(t),c(t)+\eta)$. Therefore, $u_{\delta,\varepsilon}\uparrow u_{\delta}$
as $\varepsilon\downarrow0$. 

Assume that $\Phi_{+}'(t,u_{\delta}(t))-\Phi_{-}'(t,u_{\delta}(t))<\delta$
for all $t\in E$, for some measurable set $E\subseteq H$, with non-zero
measure $\mu(E)>0$. Then we can find $\varepsilon>0$ and measurable
set $F\subseteq E$, with non-zero measure $\mu(F)>0$, such that
$\Phi_{-}'(t,x)-\Phi_{-}'(t,u_{\delta}(t))<\delta$ for all $x\in(u_{\delta}(t),u_{\delta}(t)+\varepsilon)$
and all $t\in F$. Fixed any $t\in F$ and $x\in(u_{\delta}(t),u_{\delta}(t)+\varepsilon)$,
select some $y\in(x,u_{\delta}(t)+\varepsilon)$. Then we can write
that $\Phi_{+}'(t,x)-\Phi_{-}'(t,x)\leq\Phi_{-}'(t,y)-\Phi_{-}'(t,u_{\delta}(t))<\delta$.
This result is a contradiction to the definition of $u_{\delta}$.
Consequently, $\Phi_{+}'(t,u_{\delta}(t))-\Phi_{-}'(t,u_{\delta}(t))\geq\delta$
for $\mu$-a.e.\ $t\in H$.
\end{proof}

Finally, we can prove the result stated in the beginning of this section.

\begin{proof}[Proof of \prettyref{prop:smoothness}.]
 Assume that (a)--(c) are satisfied. Let $u\in L_{0}^{\Phi}\setminus\{0\}$.
For any $k>0$ and $\lambda>1$, we can write 
\begin{align*}
\Phi^{*}(t,\Phi_{+}'(t,ku)) & \leq\Phi(t,ku)+\Phi^{*}(t,\Phi_{+}'(t,ku))=ku\Phi_{+}'(t,ku)\\
 & \leq\frac{1}{\lambda-1}\int_{ku}^{\lambda ku}\Phi_{+}'(t,x)dx\leq\frac{1}{\lambda-1}\Phi(t,\lambda ku).
\end{align*}
From (b), we obtain that $I_{\Phi^{*}}(\Phi_{+}'(t,|ku(t)|))<\infty$
for any $k>0$. Since $\Phi_{+}'(t,0)=0$ and $\Phi^{*}(t,b_{\Phi^{*}}(t))=\infty$
for $\mu$-a.e.\ $t\in T$, it follows that $I_{\Phi^{*}}(\Phi_{+}'(t,|ku(t)|))\downarrow0$
as $k\downarrow0$, and $I_{\Phi^{*}}(\Phi_{+}'(t,|ku(t)|))\uparrow\infty$
as $k\uparrow\infty$. By the continuity of $\Phi_{+}'(t,\cdot)$,
there exists only one measurable function $v$ satisfying $I_{\Phi^{*}}(v)=1$
and such that $\operatorname{sgn}v(t)=\operatorname{sgn}u(t)$ and
$|v(t)|\in\partial\Phi(t,|ku(t)|)$ for $\mu$-a.e.\ $t\in T$, for
all $k\in K(u)$. According to \prettyref{prop:smooth_point_Ku_neq_empty},
the function $u$ is a smooth point.

Conversely, let $L^{\Phi}$ be a smooth Musielak--Orlicz space. If
$E=\{t\in T:\Phi^{*}(t,b_{\Phi^{*}}(t))<\infty\}$ has non-zero measure,
then we can find a measurable set $F\subseteq E$, with $\chi_{F}\in\tilde{L}^{\Phi}\setminus\{0\}$,
and such that $I_{\Phi^{*}}(b_{\Phi^{*}}\chi_{F})<1$ and $\mu(T\setminus F)>0$.
In view of \prettyref{prop:smooth_point_Ku_eq_empty}, the function
$\chi_{F}$ is not a smooth point. This result shows that $\Phi^{*}(t,b_{\Phi^{*}}(t))=\infty$
for $\mu$-a.e.\ $t\in T$. Assume that $\Phi$ does not satisfy
the $\Delta_{2}$-condition. According to \prettyref{prop:IPhiStarLessOne},
there exists a function $u\in\tilde{L}^{\Phi}$ such that $I_{\Phi}(\lambda u)=\infty$
for any $\lambda>1$, and $I_{\Phi^{*}}(\Phi_{+}'(t,|u(t)|))<1$.
Clearly, $K(u)=\{1\}$. By \prettyref{prop:smooth_point_Ku_neq_empty},
we obtain that $u$ is not a smooth point. Thus $\Phi$ satisfies
the $\Delta_{2}$-condition.

Now suppose that $\Phi(t,\cdot)$ is not continuously differentiable
for $\mu$-a.e.\ $t\in T$. According to \prettyref{lem:u_delta},
for any $\delta>0$, the function 
\[
u_{\delta}(t)=\sup\{u\geq0:\Phi_{+}'(t,x)-\Phi_{-}'(t,x)<\delta\text{ for all }0\leq x\leq u\}
\]
is measurable. From the assumption that $\Phi(t,\cdot)$ is not continuously
differentiable for $\mu$-a.e.\ $t\in T$, we can find some $\delta_{0}>0$
for which the measurable set $H=\{t\in T:u_{\delta_{0}}(t)<\infty\}$
has non-zero measure. Denote $u=u_{\delta_{0}}$. In view of \prettyref{lem:u_delta},
we have that $\Phi_{+}'(t,u(t))-\Phi_{-}'(t,u(t))\geq\delta$, for
$\mu$-a.e.\ $t\in H$. Let $A\subseteq H$ be a measurable set,
with non-zero measure, such that $T\setminus A$ has non-zero measure,
and $I_{\Phi}(u\chi_{A})<\infty$ and $I_{\Phi^{*}}(\Phi_{+}'(t,u(t)\chi_{A}(t)))\leq1$.
We take disjoint, measurable sets $E$ and $F$, with non-zero measure,
satisfying $A=E\cup F$ and 
\[
\int_{E}[\Phi^{*}(\Phi_{+}'(t,u(t)))-\Phi^{*}(\Phi_{-}'(t,u(t)))]d\mu=\int_{F}[\Phi^{*}(\Phi_{+}'(t,u(t)))-\Phi^{*}(\Phi_{-}'(t,u(t)))]d\mu,
\]
Thus we can write 
\begin{multline*}
I_{\Phi^{*}}(\Phi_{+}'(t,u(t))\chi_{E}(t))+I_{\Phi^{*}}(\Phi_{-}'(t,u(t))\chi_{F}(t))\\
=I_{\Phi^{*}}(\Phi_{+}'(t,u(t))\chi_{E}(t))+I_{\Phi^{*}}(\Phi_{-}'(t,u(t))\chi_{F}(t))=c\leq1.
\end{multline*}
Since the set $T\setminus A$ has non-zero measure, and $\Phi^{*}(t,b_{\Phi^{*}}(t))=\infty$
for $\mu$-a.e.\ $t\in T$, we can find a sufficiently large $n_{0}\geq1$
for which $I_{\Phi^{*}}(\Phi_{+}'(t,n_{0})\chi_{T\setminus A}(t))\geq1-c$.
Let $B\subseteq T\setminus A$ be a measurable set for which $I_{\Phi^{*}}(\Phi_{+}'(t,n_{0}\chi_{B}))=1-c$.
Define the functions
\begin{align*}
v_{1} & =\Phi_{+}'(t,u(t))\chi_{E}(t)+\Phi_{-}'(t,u(t))\chi_{F}(t)+\Phi_{+}'(t,n_{0})\chi_{B}(t),\\
v_{2} & =\Phi_{-}'(t,u(t))\chi_{E}(t)+\Phi_{+}'(t,u(t))\chi_{F}(t)+\Phi_{+}'(t,n_{0})\chi_{B}(t).
\end{align*}
By $I_{\Phi^{*}}(v_{1})=I_{\Phi^{*}}(v_{2})=1$, it follows that $\Vert v_{1}\Vert_{\Phi^{*}}=\Vert v_{2}\Vert_{\Phi^{*}}=1$.
Now define $\widetilde{u}=u\chi_{A}+n_{0}\chi_{B}$. Clearly, $|f_{v_{i}}(\widetilde{u})|\leq\Vert v_{i}\Vert_{\Phi^{*}}\Vert\widetilde{u}\Vert_{\Phi,0}=\Vert\widetilde{u}\Vert_{\Phi,0}$.
In addition, 
\[
\Vert\widetilde{u}\Vert_{\Phi,0}\leq I_{\Phi}(\widetilde{u})+1=I_{\Phi}(\widetilde{u})+I_{\Phi^{*}}(v_{i})=\int_{T}\widetilde{u}v_{i}d\mu=f_{v_{i}}(\widetilde{u}),
\]
which implies that $f_{v_{i}}(\widetilde{u})=\Vert\widetilde{u}\Vert_{\Phi,0}$.
Thus $f_{v_{1}}$ and $f_{v_{2}}$ are different support functionals
at $\widetilde{u}$. This contradiction shows that $\Phi(t,\cdot)$
is continuously differentiable for $\mu$-a.e.\ $t\in T$.
\end{proof}

\bibliographystyle{plain}
\bibliography{refs_musielak-orlicz}

\end{document}